\documentclass[12pt,thmsa]{article}
\usepackage{amssymb}

\newlength{\defaultparindent}
\input{tcilatex}
\begin{document}

\author{Jean PETITOT\thanks{%
EHESS \& CREA (Ecole Polytechnique), France.\ \texttt{%
petitot@ehess.fr}.}}
\title{A transcendental view on the continuum: Woodin's conditional platonism}
\date{December 2007}
\maketitle

\noindent To appear in \emph{Intellectica}.\bigskip

\abstract{One of the main difficulty concerning the nature of the
continuum is to do justice, inside the set theoretical Cantorian framework,
to the classical conception (from Aristotle to Thom, via Kant, Peirce,
Brentano, Husserl and Weyl) according to which the continuum is a
non-compositional, cohesive, primitive, and intuitive datum. This paper
investigates such possibilities, from G\"{o}del to Woodin, of modelling inside
a \emph{ZFC}-universe the transcendence of the intuitive continuum w.r.t. its
symbolic determination.

\bigskip

}

\section{Introduction}

One of the main general philosophical problem raised by the nature of the
continuum is the conflict between two traditions: the older one, which can
be called ``neo-Aristotelian'', even if this sounds rather vague, and the
now the classical one, namely the Cantorian tradition. According to the
neo-Aristotelian tradition, the continuum is experienced and thought of as a
non-compositional, cohesive, primitive, and intuitive datum.\ It can be
segmented into parts but these parts are themselves continua and points are
only their boundaries.

This point of view was very well defended by Kant.\ As soon as in his 1770 
\textit{Dissertatio}, he emphasized the fact that

\begin{quotation}
\noindent ``a magnitude is continuous when it is not composed out of simple
elements'' (AK, II, p.~399~\footnote{%
AK refers to the collected works \emph{Kants gesammelte Schriften},
Preussische Akademie der Wissenschaften.}),
\end{quotation}

\noindent and explained that for the continuous ``pure intuitions'' of space
and time

\begin{quotation}
\noindent ``any part of time is still a time, and the simple elements which
are in time, namely the \emph{moments}, are not parts but \emph{limits}
between which a time takes place'' (AK, II, p.~399).

\noindent ``space must necessarily be conceived of as a continuous
magnitude, (...) and therefore simple elements in space are not parts but
limits'' (AK, II, p.~404).
\end{quotation}

\noindent He wrote also in the ``Anticipations of perception'' of the \emph{%
Critic of Pure Reason}:

\begin{quotation}
\noindent ``Space and time are quanta continua because no part of them can
be given without being enclosed into limits (points or moments) (...). Space
is made up only of spaces and time of times.\ Points and moments are only
limits, that is to say, simple places bounding space and time (...), and
neither space nor time can be made up of simple places, that is of integral
parts which would be given before space and time themselves'' (AK, III,
p.~154).
\end{quotation}

We see that this conception of the continuum is based on \emph{mereological}
properties, and especially on the concept of boundary (\emph{Grenze}):
points are boundaries and boundaries are dependent entities which cannot
exist independently of the entities they bound.When there are no explicit
boundaries, the continuum is characterized by the ``fusion'' of its parts.%
\footnote{%
It is very easy to construct a model of such a mereology.\ Let us take $\mathbb{%
R}$ with its standard topology and posit that the only admissible parts $U$
of $\mathbb{R}$ are its open subsets $U\in \mathcal{P}_{\func{ad}}(\mathbb{R})$.\
Let $U,V\in \mathcal{P}_{\func{ad}}(\mathbb{R})$. The complement $\lnot U$ of $%
U $ is the interior $\limfunc{Int}(\overline{U})$ of its classical
complement $\overline{U}$ (which is a closed subset).\ Therefore $U\cap
\lnot U=\emptyset $ but $U\cup \lnot U\neq \mathbb{R}$. Conversely; if $U\cup V=%
\mathbb{R}$ then $U\cap V\neq \emptyset $. This means that, for the open
mereology, $\mathbb{R}$ is undecomposable. Topological boundaries $\partial U$
are not admissible parts but only limits, and bound both $U$ and $\lnot U$.
Heyting used this mereotopology for defining truth in intuitionist logic.}
Moreover, Leibniz's principle of continuity holds (every function is
continuous).

Even after the Cantor-Dedekind arithmetization, the ``neo-Aristotelian''
non-compositionality of the continuum kept on raising fondamental problems
for some of the greatest philosophers, mathematicians, and psychologists,
such as Peirce, Brentano, Stumpf (who elaborated the key concept of
``fusion'': \emph{\ Verschmelzung}), Husserl~\footnote{%
Concerning the concept of \emph{Verschmelzung} in Stumpf and Husserl, see.
Petitot [1994].}, or Thom (for whom the continuum possessed an ontological
primacy as a qualitative homogeneous Aristotelian ``homeomer''). G\"{o}del
himself considered the real intuition of the continuum in this way and
opposed it to its set theoretical idealisation.

Peirce developed a ``synechology'', ``syneche'' being the greek term for
``continuum''. Mathematically, he was also the first, as far as I know, to
reject the continuum hypothesis CH ($2^{\aleph _{0}}=\aleph _{1}$) and to
define the power $\frak{c}$ of the continuum as a \emph{large cardinal},
namely an inaccessible cardinal (if $\kappa <\frak{c}$ then $2^{\kappa }<%
\frak{c}$). In some texts, Peirce even explained that the continuum could be
so huge that it would fail to be a cardinal.

Husserl was the first, in the third \emph{Logical Investigation}, to
formalize the idea of a mereology, and, after him, Stanis\l aw Lesniewski
developped the theory between 1916 and 1921. But the definition of
boundaries in a mereo(topo)logy remained up to now highly problematic as it
is argued, e.g., in Breysse-De Glas (2007).

All these conceptions develop the same criticism against the idea of\emph{\
arithmetizing} the continuum.\ According to them, the continuum can be
measured using systems of numbers, but no system can exhaust the substratum
it measures. They consider:

\begin{enumerate}
\item  that a point in the continuum is a discontinuity (a mark, a local
heterogeneity, a boundary) which is like a singular individuated ``atom''
which can be refered to by a symbol;

\item  that quantified sentences of an appropriate predicate calculus can be
therefore interpreted in the continuum;

\item  that systems of numbers can of course measure such systems of marks
and enable their axiomatic control;

\item  but that the arithmetization of the continuum postulates, what's
more, that the intuitive phenomenological continuum is reducible to such a
set-theoretic system (Cantor-Dedekind);

\item  and therefore that such a reductive arithmetization is unacceptable
for it violates the original intuitive mode of givenness of the continuum.
\end{enumerate}

Let us leave phenomenology and psychology for mathematics.\ Even if we adopt
a set-theoretic perspective making the continuum a set, non-compositionality
and cohesivity remain meaningful. They now mean in particular that the
continuum cannot be identified with a set of well \emph{individuated }
points. It is the case in intuitionistic logic where the law of the excluded
middle, which implies that two elements $a$ and $b$ of $\mathbb{R}$ are
different or equal, and the law of comparability, which implies $a=b$, or $%
a<b$, or $a>b$, are no longer valid. For Hermann Weyl (1918), this intrinsic
lack of individuation and localization of points in $\mathbb{R}$ characterizes
the continuum as an intuitive datum.

Given the close link discovered by Bill Lawvere between intuitionistic logic
and topos theory, it is not surprising that in many topo\"{i} the object $%
\mathbb{R}$ is undecomposable and all the morphisms $f:\mathbb{R}\rightarrow \mathbb{R%
}$ are continuous (Leibniz's principle, see e.g. John Bell's works~\footnote{%
See in particular Bell's contribution to this volume.}).

But, even in the realm of \emph{classical} logic and classical set theory
(that is $ZFC$: Zermelo-Fraenkel~+ axiom of choice) there exist many
evidences of the transcendence of the continuum relatively to its symbolic
logical control.\ Of course, they occur in non-standard models of $\mathbb{R}$,
but these non-Archimedian models are not well-founded and will not be
analyzed here.\footnote{%
For an introduction to non-standard analysis, see Petitot [1979], [1989] and
their bibliographies.} I will rather focus on their status in well-founded
models $V$ of $ZFC$.

So, the context of this paper will be the theory of the continuum in
universe of sets satisfying $ZFC$. We will see that, even in this purely
formal framework, the Kantian opposition between ``intuitive'' and
``conceptual'' remains operating, where ``conceptual'' now refers to the
logical control and symbolic determination of the continuum (as it was
already the case for Kant himself who considered arithmetic, contrary to
geometry, as ``conceptual'' and ``intellectual'', and algebra as a calculus
on ``symbolic constructions'').\footnote{\emph{Critic of Pure Reason},
``Transcendental Methodology'', AK, III, p.~471.} The key question remains
the same: is it possible to determine completely the intuitive continuum
using logical symbolic constructions?

We will meet in the sequel an alternative opposing two different types of
philosophies.

\begin{enumerate}
\item  Philosophies of the first type are ``minimizing'', ``ontologically''
deflationist (in the sense of restricting what axioms of existence are
admissible), nominalist, and constructive.\ They consider that the only
meaningful content of the continuum is the part which can be
``conceptually'' (that is symbolically) well determined and that the rest is
``inherently vague''. A celebrated representative of such a perspective is
Solomon Feferman who considers that the continuum cannot be a definite
mathematical object since some of its properties, such as Cantor's continuum
hypothesis, are not expressible by definite propositions. But, as explained
by John Steel (2004) in a criticism of Feferman:
\end{enumerate}

\begin{quotation}
\noindent ``Taken seriously, this analysis leads us into a retreat to some
much weaker constructivist language, a retreat which would toss out good
mathematics in order to save inherently vague philosophy.''
\end{quotation}

\begin{enumerate}
\item[2.]  It is why, philosophies of the second type are, on the contrary,
``maximizing'', ``ontologically'' inflationist, platonist in a sophisticated
sense, and highly non constructive.\ They aim at modelling \emph{inside} a $%
ZFC$-universe the transcendence of the intuitive continuum w.r.t. its
logical symbolic mastery.\ Owing to this, they must introduce non
constructive axioms for higher infinite.
\end{enumerate}

\section{Preliminaries\protect\footnote{%
See Jech [1978].}}

\subsection{Axioms}

We work in $ZFC$ with the classical axioms:

\emph{Extensionality}: sets are determined by their elements,

\emph{Pairing}: the pair $\left\{ a,b\right\} $ exists for every sets $a$
and $b$,

\emph{Union}: the union $\bigcup X=\left\{ u\in x\in X\right\} $ of every
set $X$ exists,

\emph{Power set}: the set $\mathcal{P}(X)=\left\{ u\subseteq X\right\} $) of
subsets of every set $X$ exists,

\emph{Comprehension} or \emph{Separation} (axiom schema): if $\varphi (x)$
is a formula, the subset $\left\{ x\in X:\varphi (x)\right\} $ exists for
every set $X$,

\emph{Replacement} (axiom schema): if $y=f(x)$ is function (i.e. a relation $%
\varphi (x,y)$ s.t. $\varphi (x,y)$ and $\varphi (x,z)$ imply $y=z$), then
the image $\left\{ f(x):x\in X\right\} $ exists for every set $X$,

\emph{Infinity}: there exists an inductive set, that is a set $I$ s.t. $%
\emptyset \in I$ and if $x\in I$ then $x\cup \{x\}\in I$,

\emph{Regularity}: all sets have minimal $\in $-elements,

\emph{Choice}: every family of sets $X_{s}$, $s\in S$, has a choice function 
$f:S\rightarrow \mathcal{P}(X)$ s.t. $f(s)\in X_{s}$ for every $s\in S$
(this axiom of existence doesn't define any specific $f$ and is highly non
constructive).

\subsection{Ordinals and cardinals}

\emph{Ordinals} are the sets that are $\in $-transitive ($y\in X$ implies $%
y\subset X$ , that is $\bigcup X\subseteq X$ or $X\subset \mathcal{P}(X)$)
and well-ordered by $\in $. All well-ordered sets are order-isomorphic to an
ordinal.\ Every ordinal is a successor: $\alpha =\beta +1$ or a limit
ordinal $\alpha =\func{Sup}\left( \beta :\beta <\alpha \right) $ (and then $%
\forall \beta <\alpha $, $\beta +1<\alpha $). A limit ordinal is like an
``horizon'' for enumeration: it is impossible to reach its limit in a finite
number of steps. The smallest limit ordinal is $\omega =\mathbb{N}$. The sum of
two ordinals is their concatenation (non commutativity: $1+\omega =\omega
\neq \omega +1$). The product of two ordinals $\alpha .\beta $ is $\beta $%
-times the concatenation of $\alpha $ (lexicographic order) (non
commutativity: $2.\omega =\omega \neq \omega .2=\omega +\omega $). An
ordinal $\alpha $ is a limit ordinal iff there exists $\beta $ s.t. $\alpha
=\omega .\beta $.

\emph{Cardinals} $\left| X\right| $ are equivalence classes of the
equivalence relation of equipotence: $X\func{eq}Y$ if there exists a
bijective (i.e. one-to-one onto) map $f:X\rightarrow Y$. They highly depend
upon the functions existing in the $ZFC$-universe under consideration.

\textbf{Cantor theorem.}\ $\left| X\right| <\left| \mathcal{P}\left(
X\right) \right| $.$\hfill \square $

Indeed, let $f:X\rightarrow \mathcal{P}\left( X\right) $.\ Then $Y=\left\{
x\in X:x\notin f(x)\right\} $ exists (Comprehension axiom). But $Y\notin
f(X) $ for if it would exist $z\in X$ with $f(z)=Y$, then $z\in
Y\Leftrightarrow z\notin Y$.\ Contradiction.$\hfill \square $

If $\left| A\right| =\kappa $, then $\left| \mathcal{P}\left( A\right)
\right| =2^{\kappa }$ since $X\subseteq A$ is equivalent to its
characteristic function $\chi _{X}:A\rightarrow \{0,1\}$ and $\chi _{X}\in
2^{\kappa }$. Cantor theorem implies therefore $\kappa <2^{\kappa }$.

For ordinals, cardinals \emph{numbers} are the $\alpha $ s.t. $\left| \beta
\right| <\left| \alpha \right| $ for every $\beta <\alpha $. They are the
minimal elements in the equivalence classes of equipotent ordinals.\ So every%
\emph{\ infinite} cardinal is a \emph{limit} ordinal. Every well-ordered set
has such a cardinal number for cardinal. These cardinal numbers are the\emph{%
\ alephs} $\aleph _{\alpha }$. Each $\aleph _{\alpha }$ has a \emph{successor%
}, namely $\aleph _{\alpha }^{+}=\aleph _{\alpha +1}$.\ If $\alpha $ is a
limit ordinal, $\aleph _{\alpha }=\omega _{\alpha }=\func{Sup}_{\beta
<\alpha }\left( \omega _{\beta }\right) $.

\textbf{Theorem.}\ \textsl{For alephs, the sum and product operations are
trivial: }$\aleph _{\alpha }+\aleph _{\beta }=\aleph _{\alpha }.\aleph
_{\beta }=\max \left( \aleph _{\alpha },\aleph _{\beta }\right) $ (the
bigger takes all)\textsl{.}$\hfill \square $

A consequence is that if $\alpha \leq \beta $ then $\aleph _{\alpha
}^{\aleph _{\beta }}=2^{\aleph _{\beta }}$.\ Indeed, $2^{\aleph _{\beta
}}\leq \aleph _{\alpha }^{\aleph _{\beta }}\leq \left( 2^{\aleph _{\alpha
}}\right) ^{\aleph _{\beta }}=2^{\aleph _{\alpha }.\aleph _{\beta
}}=2^{\aleph _{\beta }}$. But other \emph{exponentiations} raise fundamental
problems.

Under the axiom of choice $AC$, \emph{every} set can be well-ordered and
therefore all cardinals are alephs. In particular $2^{\aleph _{0}}$ is an
aleph $\aleph _{\alpha }$.

\section{The underdetermination of cardinal arithmetic in $ZFC$}

Let $V$ be a universe of set theory (i.e. a model of $ZFC$).\ We work in $%
\mathbb{R}$ or in the isomorphic Baire space $\mathcal{N}=\omega ^{\omega }$.
The first limit we meet is that the axioms of $ZFC$ are radically
insufficient for determining the cardinal arithmetic of $V$ as it is clearly
shown by the following celebrated result.\footnote{%
See Jacques Stern [1976] for a presentation.}

\subsection{Easton theorem}

For every ordinal $\alpha $ let $F(\alpha )$ be the power function defined
by $2^{\aleph _{\alpha }}=\aleph _{F(\alpha )}$.\footnote{%
More generally one can consider the power function $(\lambda ,\kappa
)\mapsto \lambda ^{\kappa }$ for each pair $(\lambda ,\kappa )$ of cardinals.%
} One can show that:

\begin{enumerate}
\item  $F$ is a monotone increasing function: if $\alpha \leq \beta $ then $%
F(\alpha )\leq F(\beta )$;

\item  \emph{K\"{o}nig's law}: $\func{cf}\left( \aleph _{F(\alpha )}\right)
>\aleph _{\alpha }$, where the \emph{cofinality} $\func{cf}(\alpha )$ of an
ordinal $\alpha $ is defined as the smallest cardinality $\chi $ of a
cofinal (i.e. unbounded) subset $X$ of $\alpha $ (i.e. $\func{Sup}X=\alpha $
). For instance, $\func{cf}\left( \omega +\omega \right) =\func{cf}\left(
\aleph _{\alpha +\omega }\right) =\omega $. Of course, $\func{cf}(\alpha )$
is a limit ordinal and $\func{cf}(\alpha )\leq \alpha $.\ The cardinal $%
\kappa $ is called \emph{regular} if $\func{cf}(\kappa )=\kappa $ i.e. if,
as far as we start with $\alpha <\kappa $, it is impossible to reach the
horizon of $\kappa $ in less than $\kappa $ steps. In some sense the length
of $\kappa $ is equal to its ``asymptotic'' length and cannot be exhausted
before reaching the horizon. As $\func{ cf}(\func{cf}(\kappa ))=\func{cf}%
(\kappa )$, $\func{cf}(\kappa )$ is always regular. As $\func{cf}\left(
\aleph _{\alpha +\omega }\right) =\omega $, $\aleph _{\alpha +\omega }$ is
always singular.\footnote{%
So a cardinal $\kappa $ is singular iff $\kappa =\bigcup_{i\in I}\alpha _{i}$
with $\left| I\right| <\kappa $ and $\left| \alpha _{i}\right| <\kappa $ $%
\forall i\in I$.}
\end{enumerate}

K\"{o}nig's law is a consequence of a generalization of Cantor theorem which
says that $1+1+\ldots $ ($\kappa $ times) $<2.2.\ldots $ ($\kappa $ times):
if $\kappa _{i}<\lambda _{i}$ $\forall i\in I$, then $\sum_{i\in I}\kappa
_{i}<\prod_{i\in I}\lambda _{i}$. Let $\kappa _{i}<2^{\aleph _{\alpha }}$
and $\lambda _{i}=2^{\aleph _{\alpha }}$ for $i<\omega _{\alpha }$. Then $%
\sum_{i<\omega _{\alpha }}\kappa _{i}<\prod_{i<\omega _{\alpha }}\lambda
_{i}=\left( 2^{\aleph _{\alpha }}\right) ^{\aleph _{\alpha }}=2^{\aleph
_{\alpha }}$.\ Therefore with an $\omega _{\alpha }$ sequence of $\kappa
_{i}<2^{\aleph _{\alpha }}$ it is impossible to get $\func{Sup}\left( \kappa
_{i}\right) =2^{\aleph _{\alpha }}$ since for infinite cardinals $\func{Sup}%
=\sum $.

An immediate corollary of K\"{o}nig's law is

\textbf{Theorem.}\ \textsl{If }$\kappa $\textsl{\ is an infinite cardinal, }$%
\kappa <\kappa ^{\func{cf}(\kappa )}$\textsl{\ (compare with Cantor: }$%
\kappa <2^{\kappa }$\textsl{).}$\hfill \square $

Indeed, if $\kappa _{i}<\kappa $ for $i<\func{cf}(\kappa )$ and $\kappa =%
\func{Sup}\left( \kappa _{i}\right) =\sum_{i<\func{cf}(\kappa )}\kappa _{i}$
then $\kappa =\sum_{i<\func{cf}(\kappa )}\kappa _{i}<\prod_{i<\func{cf}%
(\kappa )}\kappa =\kappa ^{\func{cf}(\kappa )}$.$\hfill \square $

In fact one can prove that the essential cardinals for cardinal arithmetic
are the $2^{\kappa }$ and the $\kappa ^{\func{cf}(\kappa )}$ (Gimel
function).\ They enable to compute all the $\aleph _{\alpha }^{\aleph
_{\beta }}$:

\textbf{Theorem.\footnote{%
See Jech [1978], p.\ 49.}}

\begin{enumerate}
\item  \textsl{If }$\alpha \leq \beta $\textsl{, then }$\aleph _{\alpha
}^{\aleph _{\beta }}=2^{\aleph _{\beta }}$\textsl{.}

\item  \textsl{If }$\alpha >\beta $\textsl{\ and }$\exists \gamma <\alpha $%
\textsl{\ s.t. }$\aleph _{\gamma }^{\aleph _{\beta }}\geq \aleph _{\alpha }$%
\textsl{, then }$\aleph _{\alpha }^{\aleph _{\beta }}=\aleph _{\gamma
}^{\aleph _{\beta }}$\textsl{.}

\item  \textsl{If }$\alpha >\beta $\textsl{\ and }$\forall \gamma <\alpha $%
\textsl{\ we have }$\aleph _{\gamma }^{\aleph _{\beta }}<\aleph _{\alpha }$%
\textsl{\ then}

\begin{enumerate}
\item[(a)]  \textsl{if }$\aleph _{\alpha }$\textsl{\ is regular or }$\func{cf%
}(\aleph _{\alpha })>\aleph _{\beta }$\textsl{\ then }$\aleph _{\alpha
}^{\aleph _{\beta }}=\aleph _{\alpha }$\textsl{;}

\item[(b)]  \textsl{if }$\func{cf}(\aleph _{\alpha })\leq \aleph _{\beta
}<\aleph _{\alpha }$\textsl{\ then }$\aleph _{\alpha }^{\aleph _{\beta
}}=\aleph _{\alpha }^{\func{cf}(\aleph _{\alpha })}$\textsl{.\hfill }$%
\square $
\end{enumerate}
\end{enumerate}

If the generalized continuum hypothesis ($GCH$) holds, K\"{o}nig's law is
trivial because $F(\alpha )=\alpha +1$, every cardinal $\aleph _{\alpha +1}$
is regular and therefore $\func{cf}\left( \aleph _{\alpha +1}\right) =\aleph
_{\alpha +1}>\aleph _{\alpha }$.

The fact that $ZFC$ radically underdetermines cardinal arithmetic is
particularly evident in Easton's striking result:

\textbf{Easton theorem.} \textsl{For }regular\textsl{\ cardinals }$\aleph
_{\alpha }$\textsl{, one can impose via forcing in }$ZFC$\textsl{\ the power
function }$2^{\aleph _{\alpha }}=\aleph _{F(\alpha )}$\textsl{\ for quite
every function }$F$\textsl{\ satisfying (i) and (ii).\hfill }$\square $

For regular cardinals $\kappa $, we have $\kappa ^{\func{cf}\left( \kappa
\right) }=\kappa ^{\kappa }=2^{\kappa }$ but for singular cardinals $\sigma $
we have $\sigma ^{\func{cf}\left( \sigma \right) }=\left( 
\begin{array}{c}
\sigma \\ 
\func{cf}\left( \sigma \right)
\end{array}
\right) 2^{\func{cf}\left( \sigma \right) }$ where $\left( 
\begin{array}{c}
\kappa \\ 
\lambda
\end{array}
\right) $ for $\kappa >\lambda $ is a sophisticated generalization of the
binomial formula.

The proof of Easton theorem uses iterated Cohen forcing.

\subsection{Cohen forcing}

Cohen forcing (1963)\footnote{%
Awarded a Fields Medal in 1966, Paul Cohen died on March 23, 2007.} allows
to construct in a very systematic way ``generic'' extensions $N$ of \emph{%
inner models} $M$ of $ZF$ or $ZFC$ (that is transitive $\in $-submodels $%
M\subset V$ of $ZF$ or $ZFC$ with $On\subseteq M$) where some desired
properties become valid.

Suppose for instance that, starting with a ground inner model $M$ of $ZFC$
in $V$, we want to construct another inner model $N$ where $\omega _{1}^{M}$
(that is the cardinal $\omega _{1}=\aleph _{1}$ \emph{as defined in} $M$) 
\emph{collapses} and becomes \emph{countable}.\ We need to have at our
disposal \emph{in} $N$ a surjection $f:\omega \rightarrow \omega _{1}^{M}$
that, by definition of $\omega _{1}^{M}$, cannot belong to $M$.\ Suppose
nevertheless that such an $f$ exists. Then for every $n$ the restriction $%
f\mid _{n}=\left( f(0),\ldots ,f(n-1)\right) $ exists and is an element of
the ground model $M$. Let us therefore consider the set $P=\left\{ p\right\} 
$ of finite sequences $p=\left( \alpha _{0},\ldots ,\alpha _{n-1}\right) $
of countable ordinals $\alpha _{i}<\omega _{1}^{M}$ of $M$. Such $p$ are
called \emph{forcing conditions} and must be interpreted as forcing $f\mid
_{n}=p$. The set $P$ exists, is well defined in $M$, and is endowed with a
natural partial order ``$q\leq p$ iff $p\subseteq q$''\footnote{%
That is $q<p$ means that $q$ forces a better approximation of $f$ than $p$.}%
. If $f$ exists, we can consider $G=\left\{ f\mid _{n}\right\} _{n\in \mathbb{N}%
}$ which is a subset of $P$ in $V$ s.t. $\cup G=f$.\ But as $f\notin M$, $G$ 
\emph{cannot be a subset of} $P$ \emph{in} $M$.

If $f$ exists, it is trivial to verify that $G$ satisfies the following
properties:

\begin{enumerate}
\item  \emph{Gluing and restriction conditions }(see topos theory): if $%
p,q\in G$, then $p$ and $q$ are initial segments of $f$ and are compatible
in the sense that $p\leq q$ or $q\leq p$ and therefore there exists a common
smaller element $r\in G$ satisfying $r\leq p$, $r\leq q$.

\item  for every $n\in \omega $, there exists $p\in G$ s.t. $n\in \func{dom}%
(p)$ (i.e. $\func{dom}(f)=\omega $).

\item  for every $M$-countable ordinal $\alpha <\omega _{1}^{M}$, there
exists $p\in G$ s.t. $\alpha \in \func{range}(p)$ (i.e. $\func{range}%
(f)=\omega _{1}^{M}$, it is the fundamental condition of surjectivity for
the collapsing of $\omega _{1}^{M}$ in $N$).
\end{enumerate}

Cohen's idea is to construct sets $G$ in $V$ satisfying these properties and
to show that extending the ground inner model $M$ by such a $G$ yields an
appropriate inner model $N=M\left[ G\right] $ which is the smaller inner
model of $V$ containing $M$ and $G$.

So, one supposes that a partially ordered set of forcing conditions $P$ is
given. A subset of conditions $D\subseteq P$ is called \emph{dense} if for
every $p\in P$ there is a smaller $d\leq p$ belonging to $D$. One then
defines \emph{generic} classes $G\subseteq P$ of conditions. $G\notin M$ is
generic over $M$ iff:

\begin{enumerate}
\item[(i)]  $p\in G$ and $p\leq q\in P$ implies $q\in G$,

\item[(ii)]  for every $p,q\in G$, there exists a common smaller $r\in G$
satisfying $r\leq p$, $r\leq q$,

\item[(iii)]  for every dense set $D$ of conditions $D\in M$ (be careful: $%
D\in M$, $D\subseteq P\in M$, $G\subseteq P$, but $G\notin M$), there exists 
$p\in D$ such that $p\in G$ (i.e. $G\cap D\neq \emptyset $).
\end{enumerate}

Properties (i) and (ii) mean that $G$ is a \emph{filter} for the order $\leq 
$. If $G$ is generic, the properties $(2)$ and $(3)$ above are automatically
satisfied since the sets of conditions $D_{n}=\left\{ p\in P:n\in \func{dom}
(p)\right\} $ for $n\in \omega $ and $E_{\alpha }=\left\{ p\in P:\alpha \in 
\func{range}(p)\right\} $ for $\alpha <\omega _{1}^{M}$ are dense: $(2)$
means $G\cap D_{n}\neq \emptyset $ and $(3)$ means $G\cap E_{\alpha }\neq
\emptyset $.

\textbf{Cohen's main theorem}.\ \textsl{There exists a} $ZFC$\textsl{-model} 
$\mathcal{A}=M[G]$ \textsl{such that} $(1)$ $M$ \textsl{is an inner model of}
$\mathcal{A}$, $(2)$ $G$ \textsl{is not a set in} $M$ \textsl{but is a set in%
} $\mathcal{A}$\textsl{,} $(3)$ \textsl{if} $\mathcal{A}^{\prime }$ \textsl{%
is another model statisfying }$(1)$\textsl{\ et }$(2)$\textsl{, then there
exists an elementary embedding} $j:\mathcal{A}\prec \mathcal{A}^{\prime }$ 
\textsl{such that }$j(\mathcal{A})$\textsl{\ is an inner model of }$\mathcal{%
A}^{\prime }$\textsl{\ and }$j|_{M}=\func{Id}(M)$\textsl{, }$(4)$\textsl{\ }$%
\mathcal{A}$\textsl{\ is essentially unique}.$\hfill \square $

If $j:\mathcal{A}\prec \mathcal{A}^{\prime }$ is an embedding of a model $%
\mathcal{A}$ in a model $\mathcal{A}^{\prime }$, \emph{elementarity} means
that $\mathcal{A}^{\prime }$ has exactly the same first-order theory as $%
\mathcal{A}$ in the language $\mathcal{L}_{\mathcal{A}}$ where there exist
names for \emph{every} element of $\mathcal{A}$ (that is for every set $x$
and first-order formula $\varphi $, $\mathcal{A}\models \varphi (x)$ iff $%
\mathcal{A}^{\prime }\models \varphi \left( j\left( x\right) \right) $. So,
in the first-order case, $\mathcal{A}^{\prime }$ adds only indiscernible
elements.\ A less constraining relation is \emph{elementary equivalence}: $%
\mathcal{A}$ and $\mathcal{A}^{\prime }$ are elementary equivalent, $%
\mathcal{A}\equiv \mathcal{A}^{\prime }$, if they have the same first-order
theory.\ In that case, elements of $\mathcal{A}$ characterized by
first-order sentences can be substituted for other elements of $\mathcal{A}%
^{\prime }$. It is no longer the case for an elementary embedding.

An essential feature of forcing extensions is that it is possible to
describe $M[G]$ using the language $\mathcal{L}_{G}$ which is the language $%
\mathcal{L}$ of $M$ extended by a new symbol constant for $G$.\ As was
emphasized by Patrick Dehornoy (2003), forcing is

\begin{quotation}
\noindent ``as a field extension whose elements are described by polynomials
defined on the ground field''.
\end{quotation}

In particular, the validity of a formula $\varphi $ in $M[G]$ can be coded
by a \emph{forcing relation} $p\Vdash \varphi $ \emph{defined in} $M$. The
definition of $p\Vdash \varphi $ is rather technical but an excellent
intuition is given by the idea of ``localizing'' truth, $p$ being
interpreted as a local domain (as an open set of some topological space),
and $p\Vdash \varphi $ meaning that $\varphi $ is ``locally true
everywhere'' on $p$.

\textbf{Forcing theorem}.\ \textsl{For every generic }$G\subseteq P$\textsl{%
\ , }$M[G]\models \varphi $\textsl{\ iff there exists a }$p\in G$\textsl{\
s.t. }$p\Vdash \varphi $\textsl{.\hfill }$\square $

Using forcing, we can add to $\mathbb{R}$ (i.e. to $\mathcal{P}(\omega )$) new
elements called \emph{generic reals}. Let $P$ be the partial order of binary
finite sequences $p=\left( p(0),\ldots ,p(n-1)\right) $.\ If $G\subseteq P$
is generic, $f=\cup G$ is a map $f:\omega \rightarrow \left\{ 0,1\right\} $
which is the characteristic function $f=1_{A}$ of a \emph{new} subset $%
A\subseteq \omega $ and $A\notin M$. Indeed, if $g:\omega \rightarrow
\left\{ 0,1\right\} $ defines a subset $B\subseteq \omega $ which belongs to 
$M,$ then the set of conditions $D_{g}=\left\{ p\in P:p\nsubseteq g\right\}
\in M$ is dense (if $p$ is any finite sequence it can be extended to a
sequence long enough to be different from $g$) and therefore $G\cap
D_{g}\neq \emptyset $.\ But this means $f\neq g$.

To prove the negation $\lnot CH$ of $CH$, one adds to $M$ a great number of
generic reals.\ More precisely, one embeds $\omega _{2}^{M}$ into $\left\{
0,1\right\} ^{\omega }$ (isomorphic to $\mathbb{R}$) using as forcing
conditions the set $P$ of finite binary sequences of $\omega _{2}^{M}\times
\omega $.\ If $G$ is generic, then $f=\cup G$ is a map $f:\omega
_{2}^{M}\times \omega \rightarrow \left\{ 0,1\right\} $, that is an $\omega
_{2}^{M}$-family $f=\left\{ f_{\alpha }\right\} _{\alpha <\omega _{2}^{M}}$
of generic reals $f_{\alpha }:\omega \rightarrow \left\{ 0,1\right\} $.\
Using density arguments one shows that $f$ yields an embedding $\omega
_{2}^{M}\hookrightarrow \left\{ 0,1\right\} ^{\omega }$ in $M\left[ G\right] 
$ and that $\omega _{2}^{M}$ \emph{doesn't collapse} in $M\left[ G\right] $
(because $P$ is $\omega $-saturated, i.e. there doesn't exist in $P$ any
infinite countable subset of incompatible elements).\ This implies
immediately $\lnot CH$.

Easton theorem is proved by \emph{iterating} such constructions and adding
to every regular $\aleph _{\alpha }$ as many new subsets as it is necessary
to have $2^{\aleph _{\alpha }}=\aleph _{F(\alpha )}$.

\subsection{Absoluteness}

Many philosophers and logicians which are ``deflationist'' regarding
mathematical ``ontology'' consider that the only sentences having a well
determined truth-value are those the truth-value of which is the same in all
models of $ZFC$, and that sentences the truth-value of which can change
depending on the chosen model are ``inherently vague''. Such an
antiplatonist conception has drastic consequences.\ Indeed, contrary to
first order arithmetic, which is $ZF$-absolute, that is invariant relative
to extensions of the universe (Sch\"{o}nfield theorem), all structures and
notions such as $\mathcal{N}$, $\mathbb{R}$, $\func{Card}(\chi )$, $%
x\rightarrow \mathcal{P}(x)$, $x\rightarrow \left| x\right| $, and second
order arithmetic, \emph{are not} \emph{ZF-absolute}. They can vary widely
from one model to another and can't have absolute truth value in $ZF$.\ This
``vagueness'' is one of the main classical arguments of antiplatonists
against non-constructive set theories.\ But, it has been emphasized by Hugh
Woodin in his 2003 paper \emph{Set theory after Russell.\ The journey back
to Eden} that vagueness is not an admissible argument against platonism and
shows only that it is necessary to \emph{classify} the different models of $%
ZF$ and $ZFC$. As he explained also in his talk at the Logic Colloquium held
in Paris in 2000 (quotation from Dehornoy, 2003, p.~23):

\begin{quotation}
\noindent ``There is a tendency to claim that the Continuum Hypothesis is
inherently vague and that this is simply the end of the story.\ But any
legitimate claim that $CH$ is inherently vague must have a mathematical
basis, at the very least a theorem or a collection of theorems.\ My own view
is that the independence of $CH$ from $ZFC$, and from $ZFC$ together with
large cardinal axioms, does not provide this basis. (...) Instead, for me,
the independence results for $CH$ simply show that $CH$ is a difficult
problem.''
\end{quotation}

In fact, the strong variability of the possible models of $ZFC$ is an
argument in favor of the \emph{irreducibility} of the continuum to a set of
points which can be ``\emph{individuated}'' by a symbolic description.\ 

To tackle this problem, we must look at two opposed strategies, both
introduced by G\"{o}del, one being ``minimalist'' (``ontologically''
deflationist) and the other ``maximalist'' (``ontologically'' inflationist),
and first introduce some classes of sets of reals.

\section{Borel and projective hierarchies}

In descriptive set theory, one works in $\mathbb{R}$ or in $\mathcal{N}=\omega
^{\omega }$ or in $\left\{ 0,1\right\} ^{\omega }$, and, more generally, on
metric, separable, complete, perfect (closed without isolated points) spaces 
$\mathcal{X}$ (Polish spaces). One considers in $\mathcal{X}$ different
``nicely'' definable classes of subsets $\Gamma $. The first is the \emph{Borel hierarchy} constructed from the open sets by iterating the operations
of complementation and of ``projection'' $\mathcal{X}\times \mathbb{\omega }%
\rightarrow \mathcal{X}$. If $P\subseteq \mathcal{X}\times \mathbb{\ \omega }$
(that is, if $P$ is a countable family of subsets $P_{n}\subseteq \mathcal{X}
$), one considers the subset of $\mathcal{X}$ defined by $\exists ^{\omega
}P=\left\{ x\in \mathcal{X}\ |\ \exists n\ P(x,n)\right\} $.\footnote{%
We identify predicates $\varphi (x)$, $P(x,n)$, etc. with their extensions.}
It is the union $\bigcup\limits_{n\in \omega }P_{n}$.

The $\Sigma _{1\ }^{0}$are the open subsets of $\mathcal{X}$, the $\Pi
_{1}^{0}=\lnot \Sigma _{1}^{0}$ are the closed subsets, $\Delta _{1}^{0}=\Pi
_{1}^{0}\cap \Sigma _{1}^{0}$ the clopen subsets, and the Borel hierarchy $B$
is defined by:

\[
\Pi _{n}^{0}=\left\{ \lnot \varphi \ |\ \varphi \in \Sigma _{n}^{0}\right\}
=\lnot \Sigma _{n}^{0},\ \Sigma _{n+1}^{0}=\exists ^{\omega }\lnot \Sigma
_{n}^{0}=\exists ^{\omega }\Pi _{n}^{0},\ \Delta _{n}^{0}=\Pi _{n}^{0}\cap
\Sigma _{n}^{0}. 
\]

\noindent It can be shown that this hierarchy is \emph{strict}:

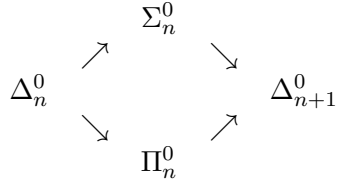
\begin{center}
\begin{tabular}{ccccc}
&  & $\Sigma _{n}^{0}$ &  &  \\ 
& $\nearrow $ &  & $\searrow $ &  \\ 
$\Delta _{n}^{0}$ &  &  &  & $\Delta _{n+1}^{0}$ \\ 
& $\searrow $ &  & $\nearrow $ &  \\ 
&  & $\Pi _{n}^{0}$ &  & 
\end{tabular}
\end{center}

One then defines the higher hierarchy of \emph{projective} sets using a
supplementary principle of construction, namely projections by continuous
projections $\mathcal{X}\times \mathcal{N}\rightarrow \mathcal{X}$, written $%
\exists ^{\mathcal{N}}$. One gets a new hierarchy beginning with the class $%
\Sigma _{1}^{1}=\exists ^{\mathcal{N}}\Pi _{1}^{0}$~-- the so called \emph{%
analytic} subsets~-- and continuing with the classes:

\[
\Pi _{n}^{1}=\left\{ \lnot \varphi \ |\ \varphi \in \Sigma _{n}^{1}\right\}
=\lnot \Sigma _{n}^{1},\ \Sigma _{n+1}^{1}=\exists ^{\mathcal{N}}\lnot
\Sigma _{n}^{1}=\exists ^{\mathcal{N}}\Pi _{n}^{1},\ \Delta _{n}^{1}=\Pi
_{n}^{1}\cap \Sigma _{n}^{1}. 
\]

\noindent For instance, $P\subseteq \mathcal{X}$ is $\Sigma _{1}^{1}\ $if
there exists a \emph{closed} subset $F\subseteq \mathcal{X}\times \mathcal{N}
$\ such that: $P(x)\Leftrightarrow \exists \alpha \ F(x,\alpha ).$ In the
same way, $P\subseteq \mathcal{X}$ is $\Sigma _{2}^{1}\ $if there exists an 
\emph{open} subset\ $G\subseteq \mathcal{X}\times \mathcal{N}\times \mathcal{%
N}$\ such that: $P(x)\Leftrightarrow \exists \alpha \ \forall \beta \
G(x,\alpha ,\beta )$,$\ $etc.

More generally, one can define projective sets in $V$ using the cumulative
hierarchy of successive levels of $V$ indexed by the class $On$ of ordinals:$%
\;V_{0}=\emptyset $, $V_{\alpha +1}=\left\{ x:x\subset V_{\alpha }\right\} $
for a successor ordinal, and $V_{\lambda }=\bigcup_{\alpha <\lambda
}V_{\alpha }$ for $\lambda $ a limit ordinal. Then $P$ is projective if it
is definable with parameters over $\left( V_{\omega +1},\in \right) $. More
precisely, $P\subset V_{\omega +1}$ is $\Sigma _{n}^{1}$ if it is the set of
sets $x$ s.t. $\left( V_{\omega +1},\in \right) \models \varphi (x)$ for a $%
\Sigma _{n}$ formula $\varphi (x)$, that is a formula of the form $\varphi
(x)=\exists x_{1}\forall x_{2}\ldots \psi $ with $n$ quantifiers and a $\psi 
$ having only \emph{bounded} quantifiers.\footnote{%
Bounded quantifiers are of the form $\exists y\in z$ and $\forall y\in z$.}

As the Borel hierarchy, the projective hierarchy is \emph{strict} and it is
a continuation of the Borel hierarchy according to:

\textbf{Suslin theorem.} $B=\Delta _{1}^{1}$\textsl{.\hfill }$\square $

This theorem can be interpreted as a \emph{construction principle}: it\emph{%
\ }asserts that the complex operation of continuous projection can be
reduced to an iteration of simpler operations of union and complementation.

There exist strict $\Pi _{n}^{1}$ and $\Sigma _{n}^{1}$ sets, which are very
natural in classical analysis. For instance, in the functional space $C[0,1]$
of real continuous functions on $[0,1]$ endowed with the topology of uniform
convergence, the subset 
\[
\left\{ f\in C\left[ 0,1\right] \ |\ f\ \text{smooth}\right\} 
\]

\noindent is $\Pi _{1}^{1}$ (but not $\Delta _{1}^{1}$). In the space $%
C[0,1]^{\omega }$ of countable sequences $(f_{i})$ of functions, the subset:

\[
\left\{ (f_{i})\in C\left[ 0,1\right] ^{\omega }\ \left| 
\begin{array}{c}
(f_{i})\ \text{converges\ for the topology } \\ 
\text{of simple convergence}
\end{array}
\right. \right\} 
\]

\noindent is $\Pi _{1}^{1}$, and the subset:

\[
\left\{ (f_{i})\in C\left[ 0,1\right] ^{\omega }\ \left| 
\begin{array}{c}
\text{a sub-sequence\ converges\ for the} \\ 
\text{topology of simple convergence}
\end{array}
\right. \right\} 
\]

\noindent is $\Sigma _{2\ }^{1}$and every $\Sigma _{2\ }^{1}$ can be
represented that way (Becker [1992]):

\textbf{Becker representation theorem.} \textsl{For every }$\Sigma _{2\
}^{1} $\textsl{-set }$S\subseteq C[0,1]$\textsl{\ there exists a sequence }$%
(f_{i}) $\textsl{\ such that } 
\[
S=A_{(f_{i})}=\left\{ g\in C\left[ 0,1\right] \ \left| 
\begin{array}{c}
\text{a sub-sequence of}\ (f_{i})\text{\ converges} \\ 
\text{towards\ }g\text{ for the topology} \\ 
\text{of simple convergence}
\end{array}
\right. \right\} . 
\]

Another examples are given by the compact subsets $K\in \mathcal{K}\left( 
\mathbb{R}^{n}\right) $ of $\mathbb{R}^{n}$: for $n\ge 3$, 
\[
\left\{ K\in \mathcal{K}\left( \mathbb{R}^{n}\right) \ |\ K\ \text{arc connected%
}\right\} 
\]

\noindent is $\Pi _{2}^{1}$, and for $n\ge 4$, 
\[
\left\{ K\in \mathcal{K}\left( \mathbb{R}^{n}\right) \ |\ K\ \text{simply\
connected}\right\} 
\]

\noindent is also $\Pi _{2}^{1}$.

In fact, projective sets can be considered as the ``reasonably'' definable
subsets of $\mathbb{R}$.

\section{The ``minimalist'' strategy of the constructible universe}

The first G\"{o}delian strategy for constraining the structure of $ZF$%
-universes consisted in \emph{restricting} the universe $V$. It is the
strategy~-- refered to as $V=L$~-- of \emph{constructible} sets (G\"{o}del
1938).

To define $L$ one substitutes, in the construction of the cumulative
hierarchy $V_{\alpha }$ of $V$ by means of a transfinite recursion on the $%
x\rightarrow \mathcal{P}(x)$ operation, the power sets $\mathcal{P}(x)$~--
which are not $ZF$-absolute~-- with \emph{smaller} sets $\mathcal{D}%
(x)=\{y\subseteq x~|$~$y$ elementary$\}$ (where ``elementary'' means
definable by a first order formula over the structure $\left\langle x,\in
,\{s|s\in x\}\right\rangle $)~-- which are $ZF$-absolute. $L$ is then
defined as $V$ using a transfinite recursion on ordinals: $L_{0}=\emptyset $
, $L_{\alpha +1}=\mathcal{D}(L_{\alpha })$, $L_{\lambda
}=\bigcup\limits_{\alpha <\lambda }L_{\alpha }$ if $\lambda $ is a limit
ordinal, and $L=\bigcup\limits_{\alpha \in On}L_{\alpha }$. The absoluteness
of $L$ comes from the fact that each level $L_{\alpha }$ is constructed
using only unambiguous formulae and parameters belonging to the previous
stages $L_{\beta }$, $\beta <\alpha $.

G\"{o}del (1938, 1940) has shown that if $V=L$ it is possible to define a 
\emph{global wellordering} on $L$, which is a very strong form of \emph{global} $AC$. The wellorder relation is defined by a transfinite induction
on the levels $\alpha $.\ If $x$ and $y$ are of different levels their order
is the order of their respective levels.\ If they are of the same level,
their order is first that of the G\"{o}del numbers of their minimal defining
formulae, and then the order of their parameters (which are of lower order
and therefore wellordered by the induction hypothesis).\ G\"{o}del also
proved that in $ZF$ we have $\left( V=L\right) \vdash GCH$.

$L$ is in fact the \emph{smallest} inner model of $V$:

\begin{enumerate}
\item[(i)]  $On\subset L$,

\item[(ii)]  $L$ is \emph{transitive}: if $y\in _{V}x$ and $x\in _{L}L$,
then $y\in _{L}L$,

\item[(iii)]  $(L,\in _{L})$ is a model of $ZF$.
\end{enumerate}

\noindent It can be defined in $V$ by a statement $L(x)=\,$``$x$ is
constructible'' which is \emph{independent} of $V$ ($ZF$-absolute), and in
that sense, it is a \emph{canonical} model of $ZFC$.

\textbf{Remark.}\ It must be emphasized that the constructible universe $L$
is not constructive since it contains the class $On$ of ordinals which is
non constructive.\ But the characteristic property of $L$ is that it reduces
non-constructivity exactly to $On$.

In the constructible universe $L$ there exists a $\Delta _{2}^{1}$-wellorder
relation $<$ on $\mathbb{R}$. According to a theorem due to Fubini, such a
wellordering \emph{cannot be Lebesgue measurable} and there exist therefore
in $L$ $\Delta _{2}^{1}$ sets which, despite the fact they belong to the low
levels of the projective hierarchy and are ``simple'' and ``nice'' to
define, are nevertheless not Lebesgue measurable and therefore not
well-behaved.

With regards to $CH$, one uses the fact that the $\Delta _{2}^{1}$-wellorder
relation $<$ on $\mathbb{R}$ is a fortiori $\Sigma _{2}^{1}$, and that the $%
\Sigma _{2}^{1}$ are the $\aleph _{1}$-Suslin sets.\ If $\chi $\ is an
infinite cardinal, $P\subseteq \mathbb{R}$\ is called a $\chi $-Suslin set if
it exists a closed subset $F\subseteq \mathbb{R}\times \chi ^{\mathbb{N}}$\ s.t. $%
P=\exists ^{\chi ^{\mathbb{N}}}F$\ (i.e. $P$\ is the projection of $F$). The $%
\Sigma _{1}^{1}$ are, by definition, the $\aleph _{0}$-Suslin sets. Indeed,
if $\chi =\aleph _{0}$ then $P=\exists ^{\mathbb{R}}F$ and therefore $P\in
\Sigma _{1}^{1}$. A theorem of Martin says that $P\subseteq \mathcal{X}$ is
an $\aleph _{n}$-Suslin set iff $P=\bigcup\limits_{\xi <\aleph _{n}}P_{\xi }$
with $P_{\xi }$ Borelians.\footnote{%
See Moschovakis [1980], p. 97.} As the wellordering $<$ on $\mathbb{R}$ is $%
\Sigma _{2}^{1}$, according to a theorem of Sch\"{o}nfield, its ordinal is $%
<\aleph _{2}$ and $CH$ is therefore valid.

In spite of its intrinsic limitations, $L$ is a very interesting model of $%
ZFC$, which possesses a ``fine structure'' interpolating between the
different $F_{\alpha }$ and very rich combinatorial properties investigated
by Jensen. One of its main properties is the following.\ Let us first define
what is a \emph{club} (``closed unbounded'' subset) $C\subseteq \alpha $ of
a limit ordinal $\alpha $: $C$ is closed for the order topology (i.e. limits
in $C$ belong to $C$: if $\beta <\alpha $ and $\func{Sup}\left( C\cap \beta
\right) =\beta $, then $\beta \in C$) and unbounded in $\alpha $ (for every $%
\beta <\alpha $) there exists an element $\gamma \in C$ s.t. $\beta <\gamma $
.\ For a cardinal $\kappa $, let $\square _{\kappa }$ be the property that
there exists a sequence of clubs $C_{\alpha }\subseteq \alpha $ with $\alpha
<\kappa ^{+}$ s.t. $C_{\alpha }$ is of order type $\leq \kappa $ (and $%
<\kappa $ if $\func{cf}\left( \alpha \right) <\kappa $) and if $\lambda $ is
a limit point of $C_{\alpha }$ then $C_{\alpha }\cap \lambda =C_{\lambda }$. 
$\square _{\kappa }$ is used to construct systematically and coherently
bijections between $\kappa $ and ordinals $\kappa \leq \alpha <\kappa ^{+}$
by cofinalizing the $\alpha $ by clubs.\ We have:

\textbf{Theorem (Jensen, 1970).}\ $V=L\models \forall \kappa \;\square
_{\kappa }$.$\hfill \square $

$\square _{\kappa }$ constrains the structure of the \emph{stationary}
subsets $S$ of $\kappa ^{+}$ ($S\subset \kappa ^{+}$ is stationary if $S\cap
C\neq \emptyset $ for every club $C$ in $\kappa ^{+}$). These cannot reflect
at some ordinal $\alpha <\kappa $ of $\func{cf}\left( \alpha \right) >\omega 
$, where ``reflect'' means ``remaining stationary in $\alpha $''.

One can generalize the concept of constructibility in two ways.\footnote{%
See Kanamori [1994], p.~34.}\ First, if $A$ is any set, one can relativize
definability to $A$ taking $\mathcal{D}_{A}(x)=\{y\subseteq x~|$~$y$
definable by a first order formula of the structure $\left\langle x,\in
,\{s|s\in A\cap x\}\right\rangle $)$\}$.\ One gets that way the universe,
called $L\left[ A\right] $, of constructible sets relative to $A$. In $%
L\left[ A\right] $ the only remaining part of $A$ is $A\cap L\left[ A\right]
\in L\left[ A\right] $. As $L$, $L\left[ A\right] $ satisfies $AC$ and is $%
ZF $-absolute.\ On the other hand, one can start the recursive construction
of $L$ not with $L_{0}=\emptyset $ but with the transitive closure of $%
\left\{ A\right\} $, $L_{0}\left( A\right) $.\ One gets that way $L\left(
A\right) $ which is the smallest inner model containing $On$ and $A$. If
there is a wellordering on $A$ (it the case if $AC$ is valid), then $L\left(
A\right) $ is globally wellordered for the same reasons as $L$. In
particular, $L(\mathbb{R})$ is a good compromise between the non
constructibility of $\mathbb{R}$ and the constructibility from $\mathbb{R}$ of the
rest of the universe.

In spite of its interest, the structure of $L$ is rather pathological with
regards to the continuum and many of the above results are in some sense
counterintuitive. They result from the fact that the $AC$, which implies the
existence of very complicated and irregular sets, remains valid in $L$ and
that the axiom of constructibility $V=L$ forces some of them to exist \emph{%
inside} the projective hierarchy which should be composed only of relatively
simple and regular sets: nicely definable sets are not necessarily
well-behaved.

It is the reason why many specialists consider that the strategy $V=L$ is
dramatically too restrictive and, moreover, that its restriction to
constructibility is not philosophically justifiable.\ For instance, John
Steel (2000) claims:

\begin{quotation}
\noindent ``The central idea of descriptive set theory is that definable
sets of reals are free from the pathologies one gets from a wellorder of the
reals.\ Since $V=L$ implies there is a $\Delta _{2}^{1}$ wellorder of the
reals, under $V=L$ this central idea collapses low in the projective
hierarchy, and after that there is, in an important sense, \emph{no}
descriptive set theory. One has instead infinitary combinatorics on $\aleph
_{1}$.\ This is certainly not the sort of theory that looks useful to
Analysts.''
\end{quotation}

One could think that generalizations of constructibility such as $L\left(
A\right) $ or $L\left[ A\right] $ would overcome the problem.\ But it is not
the case.

\section{The ``maximalist'' strategy of large cardinals}

It is therefore justified to reverse the strategy and to look for\emph{\
additional} axioms, which could be considered ``natural'', for $ZF$ and $ZFC$%
, and to try to \emph{generalize} to such augmented axiomatics the search of
canonical models and fine combinatorial structures. As was emphasized by
John Steel (2004):

\begin{quotation}
\noindent ``In extending $ZFC$, we are attempting to \emph{maximize
interpretative power}''.
\end{quotation}

\noindent And there is place for philosophy in such a maximizing strategy
since the problem is not only to find a solution to the continuum problem
but also to understand what ``to be a solution'' means. By the way, to study
such ``maximizing'' large cardinals models is perfectly compatible with a
minimalist perspective: one has only to relativize the theory to the
constructible subuniverse $L$ since $ZFC+``V=L"\vdash \varphi $ is
equivalent to $ZFC\vdash \varphi ^{L}$. As explained by Steel (2004),
suppose that the philosopher $A$ believes in $L$ and the philosopher $B$ in $%
L[G]$ with $G$ forcing the adjunction of $\omega _{2}$ reals to the model of 
$\mathbb{R}$ in $L$. $A$ believes in $CH$ and $B$ in $\lnot CH$, but $B$ can
interpret the formulae $\varphi $ of $A$ as its own $\varphi ^{L}$ and $A$
can interpret the formulae $\varphi $ of $B$ as forced $\varphi $ (the truth
of $\Vdash \varphi $ being definable in the ground model $L$, see above).
There is therefore no real conflict.

Different ``maximizing'' strategies have been considered:

\begin{enumerate}
\item  Iterate transfinitely theories $T_{\alpha +1}=T_{\alpha }+$
``consistency of $T_{\alpha }$'' starting from $ZF$ or $ZFC$.

\item  Postulate ``good'' regularity properties of projective sets, and
therefore of the continuum.

\item  Make the theory of the continuum ``\emph{rigid}'', that is define 
\emph{under which conditions the properties of} $\mathbb{R}$ \emph{cannot be
further modified by forcing}.
\end{enumerate}

Strategy (3) tries to \emph{reduce}~-- and even to neutralize~-- the
variability induced by forcing.\ The ideal aim would be forcing invariance
to make the theories of $\mathbb{R}$ and $\mathcal{P}(\mathbb{R})$ in some sense
as ``rigid'' as first order arithmetic. It is an extremely difficult program
and we will first evoke some classical results concerning $\mathbb{R}$. $CH$
concerns $\mathcal{P}(\mathbb{R})$ the forcing invariance of which is the
object of more recent works of Woodin. But we first emphasize the fact that
strategies ($1$), ($2$) converge towards the introduction of \emph{large
cardinal axioms} ($LCAs$) which express the existence of higher infinities.
Indeed, it seems that every ``maximizing'' strategy is in some sense
equivalent to a $LCA$. Look for instance at the \emph{Proper Forcing Axiom} $%
PFA$. A forcing $P$ is called \emph{proper} if, for every regular
uncountable cardinal $\lambda $, it preserves the stationary subsets of $%
\left[ \lambda \right] ^{\omega }$ (the set of countable subsets of $\lambda 
$).

\textbf{Proper Forcing Axiom.}\ If the forcing $P$ is proper and if the $%
D_{\alpha }$'s are dense subsets indexed by the countable ordinals $\alpha
<\omega _{1}$, then there exists a filter $G\subseteq P$ intersecting all
the $D_{\alpha }$'s. (Compare with the definition of $G$ being generic).$%
\hfill \square $

Many results are known for $PFA$.\ It implies $2^{\aleph _{0}}=\aleph _{2}$
(and therefore $\lnot CH$), it implies projective determinacy $PD$ (Woodin)
and $AD$ (Steel, 2007) for the inner model $L\left( \mathbb{R}\right) $.%
\footnote{%
See below \S \S\ $7.2$ and $8$ for a definition of $PD$ and $AD$.} As far as
its consistency strength is concerned, it is known that 
\[
\func{Con}\left( \exists \kappa \text{ supercompact}\right) \Rightarrow 
\func{Con}\left( PFA\right) \Rightarrow \func{Con}\left( \exists \kappa 
\text{ measurable}\right) . 
\]

\noindent It is conjectured that in fact $PFA$ is equiconsistant with ``$%
\exists \kappa $ supercompact'' (see below for a definition of
supercompacity).

It also seems that there exists a wellordering of $LCAs$ which can be
defined by the inclusion of their sets of $\Sigma _{2}^{1}$ (and even $\Pi
_{1}^{0}$) consequences. As emphasized by John Steel (2000):

\begin{quotation}
\noindent ``It seems that the consistency strengths of all natural
extensions of $ZFC$ are wellordered, and the large cardinal hierarchy
provides a sort of yardstick which enables us to compare these consistency
strengths.''
\end{quotation}

Philosophically speaking, the nominalist confusion between a strong
``quasi-ontology'' for sets and a realist ``true'' ontology of abstract
idealities has disqualified such axioms.\ But I think that such a dogmatic
prejudice has been a great philosophical mistake. Indeed, I think that one
of the best philosophical formulation of incompleteness is precisely to say
that a ``good'' theory of the continuum requires a very strong
``quasi-ontology'' for sets, a \emph{maximal one}, not a minimal one. A
``good'' regularity of the continuum entails \emph{for objective reasons} a
strong ``platonist'' commitment concerning higher infinities. This key point
has been perfectly emphasized by Patrick Dehornoy:

\begin{quotation}
\noindent ``properties which put into play objects as `small' as sets of
reals (...) are related to other properties which put into play very `huge'
objects which seem very far from them.''\footnote{%
Dehornoy [1989].}
\end{quotation}

Some specialists call ``reverse descriptive set theory'' this remarkable
equivalence between properties of regularity of projective sets and $LCAs$.

There are many theorems showing that the platonist ``cost'' of a ``good''
theory is very high. Let us for instance mention one of the first striking
theorems proved by Robert Solovay using forcing. Let $CM$ be the axiom of
existence of a measurable cardinal (see below for a definition).

\textbf{Solovay theorem (1969).} $ZFC+CM\vdash $ \textsl{every} $\Sigma
_{2}^{1}$\textsl{\ is ``regular'' (where ``regular'' means properties such
as Baire property, Lebesgue measurability, and perfect set property).}%
\footnote{%
See Moschovakis [1980], p. 284.}$\hfill \square $

\section{Regularity of projective sets}

\subsection{The regularity of analytic sets}

The French school (Borel, Baire, Lebesgue) and the Russian and Polish
schools (Suslin, Luzin, Sierpinski) initiated the study of the Borel and
projective classes and achieved deep results concerning their 
\emph{regularity} and their\emph{\ representation} where ``regularity'' means
Lebesgue measurability, or the perfect set property (to be countable or to
contain a perfect subset, i.e. a closed subset without isolated point), or
the Baire property (to be approximated by an open subset up to a meager set,
i.e. a countable union of nowhere dense sets).

The first regularity theorem is the celebrated:

\textbf{Cantor-Bendixson theorem.} \textsl{If }$A\subseteq \mathbb{R}$ \textsl{%
is }closed\textsl{, then\ }$A$\textsl{\ can be decomposed in a unique way as
a disjoint union }$A=P+S$\textsl{\ where }$P$\textsl{\ is }perfect\textsl{\
and }$S$\textsl{\ }countable\textsl{.\hfill }$\square $

As a perfect set $P$ is of cardinality $\left| P\right| =2^{\aleph _{0}}$,
the\emph{\ continuum hypothesis }$CH$\emph{\ holds for the closed sets }$\Pi
_{1}^{0}$.

Another early great classical theorem of regularity is the:

\textbf{Suslin theorem.} \textsl{The analytic subsets }$\Sigma _{1}^{1}$%
\textsl{\ shares the perfect subset property and }$CH$\textsl{\ is therefore
true for the }$\Sigma _{1}^{1}$ \textsl{sets.\hfill }$\square $

In the same way, one can show that the $\Sigma _{1}^{1}$ share the Baire
property and that the $\Sigma _{1}^{1}$ and $\Pi _{1}^{1}$ are Lebesgue
measurable. But it is \emph{impossible} to show in $ZF$ that the $\Delta
_{1}^{1}$ and $\Sigma _{2}^{1}$ share the perfect set property and to show
in $ZFC$ that the $\Delta _{2}^{1}$ share the Baire property. In fact many
of the ``natural'' properties of the projective sets go \emph{far beyond}
the demonstrative strength of $ZF$ and $ZFC$. It is therefore
methodologically and philosophically justified to look for additional axioms.

\subsection{Projective determinacy and the ``regularity'' of the continuum}

A very interesting regularity hypothesis is the so called \emph{determinacy }
property. One considers infinite games on sets $X$. Each player (I and II)
plays in turn an element $a$ of $X$:\bigskip

$
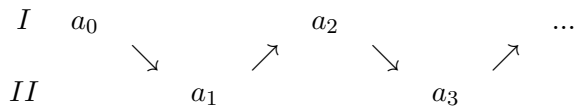
\begin{tabular}{cccccccccc}
$I$ & $a_{0}$ &  &  &  & $a_{2}$ &  &  &  & ... \\ 
&  & $\searrow $ &  & $\nearrow $ &  & $\searrow $ &  & $\nearrow $ &  \\ 
$II$ &  &  & $a_{1}$ &  &  &  & $a_{3}$ &  & 
\end{tabular}
\bigskip $

\noindent At the end of the game we get a sequence $f\in X^{\mathbb{N}}$. Let $%
A\subset X^{\mathbb{N}}$. The player I (resp. II) wins the play $f$ of the game 
$G=G_{X}(A)$ associated to $A$ if $f\in A$ (resp. if $f\notin A$).

\textbf{Definition.} $A$\textsl{\ is called }determined\textsl{\ (written }$%
\func{Det}(A)$\textsl{\ or }$\func{Det}G_{X}(A)$\textsl{) if one player has
a winning strategy. Therefore }$A$\textsl{\ is determined iff } 
\[
\exists a_{0}\forall a_{1}\exists a_{2}...(a_{0},a_{1},a_{2},...)\in A. 
\]

Determinacy is a strong property of ``regularity''. Indeed, for every $%
A\subset \mathbb{R}$ ($\mathbb{R}$ being identified with $\mathcal{N}=\omega
^{\omega }$), $\func{Det}(A)\Rightarrow $``$A$ satisfies the Baire and the
perfect subset properties, and is Lebesgue measurable''.

The first theorem linking determinacy with the projective hierarchy has been
the key result:\footnote{%
See Grigorieff [1976] and Moschovakis [1980], p. 288.}

\textbf{\ Gale-Stewart Theorem(1953).} $ZFC\vdash $\textsl{\ closed subsets }%
$A$\textsl{\ of }$X^{\mathbb{N}}$\textsl{\ (the }$\Pi _{1}^{0}$\textsl{) are
determined.\hfill }$\square $

After many efforts, Donald Martin proved a fundamental theorem which
concluded a first stage of the story:

\textbf{Martin theorem (1975).} $ZFC\vdash $\textsl{\ Borel sets (the }$%
\Delta _{1}^{1}$\textsl{) are determined.\hfill }$\square $

This celebrated result shows that $ZFC$ is a ``good'' axiomatic for the
Borel subsets of $\mathbb{R}$.\ But, it is the \emph{limit} of what is provable
in $ZFC$. Indeed, $ZFC$ \emph{cannot} imply the determinacy of $\Sigma
_{1}^{1}$-sets since in the constructible model $L$ of $ZFC$ there exist $%
\Sigma _{1}^{1}$-sets that don't share the perfect set property. As for $\Pi
_{1}^{1}$-sets, their determinacy implies the measurability of the $\Sigma
_{2}^{1}$-sets, but in $L$ there exists a $\Delta _{2}^{1}$-wellorder of $%
\mathbb{R}$, which, according to Fubini theorem, cannot be Lebesgue measurable.

\section{The necessity of large cardinals and ``reverse'' descriptive set
theory}

To prove determinacy results for projective sets beyond $\Delta _{1}^{1}$,
one must introduce additional axioms and many converging results show that
the most natural are large cardinal axioms. The first example was introduced
by Stan Ulam.\ If $X$ is a set, a \emph{filter} $\mathcal{U}$ over $X$ is a
set of subsets of $X$, $\mathcal{U}\subseteq \mathcal{P}(X)$, s.t. (i) $%
\emptyset \notin \mathcal{U}$, (ii) if $U\in \mathcal{U}$ and $U\subseteq V$
then $V\in \mathcal{U}$, (iii) if $U,V\in \mathcal{U}$ then $U\cap V\in 
\mathcal{U}$ (i.e., the complementary set of $\mathcal{U}$ in the Boolean
algebra $\mathcal{P}(X)$ is an ideal). $\mathcal{U}$ is an \emph{ultrafilter}
if it is maximal, namely if for every $U\subseteq X$, either $U\in \mathcal{U%
}$ or $X-U\in \mathcal{U}$.\ For every $x\in X$, $\mathcal{U}_{x}=\left\{
U\subseteq X:x\in U\right\} $ is an ultrafilter called ``principal''. A non
principal ultrafilter is called ``free''.

\textbf{Definition.} \textsl{A cardinal} $\chi >\omega $ \textsl{is}
measurable \textsl{if it bears a free ultrafilter }$\mathcal{U}$\textsl{\
which is }$\chi $\textsl{-}complete \textsl{(that is stable w.r.t. }$\chi $%
\textsl{-infinite intersections }$\bigcap\limits_{\lambda <\chi }X_{\lambda
} $\textsl{\ with }$\lambda <\chi $\textsl{). It is equivalent to say that }$%
\chi $\textsl{\ bears a} measure $\mu $ \textsl{\ with range }$\left\{
0,1\right\} $\textsl{\ (with }$\mu (\chi )=1$\textsl{), diffuse (without
atoms: }$\forall \xi \in \chi $ \textsl{we have} $\mu (\{\xi \})=0$\textsl{)
and }$\chi $\textsl{-additive. The equivalence is given by }$\mu
(A)=1\Leftrightarrow A\in \mathcal{U}$\textsl{\ and }$\mu
(A)=0\Leftrightarrow \chi -A\in \mathcal{U}$ \textsl{and is analogous to the
opposition between finite and infinite subsets in }$\omega $\textsl{.\hfill }%
$\square $

A first typical result was another theorem due to Donald Martin:

\textbf{Martin theorem (1970).} $ZFC+MC\vdash \func{Det}\left( \Sigma
_{1}^{1}\right) .\hfill \square $

\textbf{Corollary: Solovay theorem (1969).}\textsl{\ }$ZFC+MC\vdash $\textsl{%
\ \ ``the }$\Sigma _{2}^{1}$\textsl{-sets\ are `regular' ''.}$\hfill \square 
$

But Solovay also showed that $ZFC+MC\nvdash PD$ (where $PD$ is the
axiom of \emph{Projective Determinacy}: every projective $A\subseteq \mathbb{R}$
is determined, see below) since $ZFC+PD\vdash \func{Cons}(ZFC+MC)$ and
therefore if $ZFC+MC\vdash PD$ we would have $ZFC+MC\vdash \func{Cons}%
(ZFC+MC)$, which would contradict G\"{o}del theorem.

\textbf{Scott theorem (1961)}. $MC$ \textsl{is false} \textsl{in }$V=L$%
\textsl{\ and therefore }$ZFC\nvdash CM$\textsl{.\hfill }$\square $

The proof of Scott theorem uses the concept of an \emph{ultrapower} $V^{%
\mathcal{U}}$ where $\mathcal{U}$ is an ultrafilter on a set $S$.\ The
elements of $V^{\mathcal{U}}$ are the maps $f:S\rightarrow V$, $f$ and $g$
being equivalent if they are equal almost everywhere (a.e.), that is if $%
\left\{ s\in S:f(s)=g(s)\right\} \in \mathcal{U}$.\ Any element $x$ of $V$
is represented by the constant map $f(s)=x$ and this defines a canonical
embedding $j:V\hookrightarrow V^{\mathcal{U}}$. If $\varphi \left(
x_{1},\ldots ,x_{n}\right) $ is a formula of the language of $V$, $\varphi
\left( f_{1},\ldots ,f_{n}\right) $ is valid in $V^{\mathcal{U}}$ $\left( V^{%
\mathcal{U}}\models \varphi \left( f_{1},\ldots ,f_{n}\right) \right) $ iff $%
\varphi \left( f_{1},\ldots ,f_{n}\right) $ is valid a.e., that is if 
\[
\left\{ s\in S:V\models \varphi \left( f_{1}(s),\ldots ,f_{n}(s)\right)
\right\} \in \mathcal{U}. 
\]

One shows that $V^{\mathcal{U}}$ is well-founded if the ultrafilter $%
\mathcal{U}$ is $\omega _{1}$-complete and that there exists an isomorphism
between $\left\langle V^{\mathcal{U}},\in _{\mathcal{U}}\right\rangle $ and $%
\left\langle M_{\mathcal{U}},\in \right\rangle $ where $M_{\mathcal{U}}$ is
an inner model (Mostowski collapsing lemma).\ A fundamental theorem of \L
o\v {s} says that $j:V\prec M_{\mathcal{U}}$ is an elementary embedding (see
above for the definition). If $j$ is an elementary embedding $j:M\prec M^{*}$
of models of $ZFC$ where $M^{*}$ is an inner model of $M$ and if $\alpha \in
On(M)$ is an ordinal in $M$, one has $j(\alpha )\in On(M^{*})\subset On(M)$
and, because of the elementarity of $j$, $\alpha <\beta \Leftrightarrow
j(\alpha )<j(\beta )$. This implies $j(\alpha )\ge \alpha $. One shows that
there exists necessarily an ordinal $\alpha $ s.t. $j(\alpha )>\alpha $. Let 
$\chi $ be the smallest of these $\alpha $. It is called the \emph{critical
ordinal }$\func{crit}(j)$ of $j$.

\textbf{Theorem}.\ \textsl{If the free ultrafilter }$\mathcal{U}$\textsl{\
on the measurable cardinal }$\chi $\textsl{\ is }$\chi $\textsl{-complete,
then }$\func{crit}(j)=\chi $\textsl{\ and therefore }$j\left( \chi \right)
>\chi $\textsl{.\hfill }$\square $

\textbf{Corollary: Scott theorem.\hfill }$\square $

Indeed, suppose there exists a $MC$ and let $\chi $ be the least $MC$. Now
suppose that $V=L$. Elementarity implies $M_{\mathcal{U}}=L$ since $M_{%
\mathcal{U}}$ is an inner model satisfying the axiom of constructibility and
therefore $L\subseteq M_{\mathcal{U}}\subseteq V=L$. Then in $M_{\mathcal{U}%
}=L=V$, $j\left( \chi \right) $ is the least $MC$, which contradicts $%
j\left( \chi \right) >\chi $.$\hfill \square $

Measurable cardinals $\chi $ are very large; such a $\chi $ is \emph{regular}
(there exists no unbounded $f:\lambda \rightarrow \chi $ with $\lambda <\chi 
$), \emph{strongly inaccessible} ($\forall \lambda <\chi $, $2^{\lambda
}<\chi $), and preceded by $\chi $ strongly inaccessible smaller cardinals.
But as large as they may be, $MCs$ guarantee only the determinacy of the
lowest post-Borelian level of definable subsets of $\mathbb{R}$. To guarantee
the determinacy of all projective subsets, one needs much stronger axioms,
such as $PD$, which are not entailed by $MC$ (see Solovay's remark above)

As we will see below, many specialists consider that Projective Determinacy
is a ``good'' axiomatic for $\mathbb{R}$.$\;$Indeed, $PD$ is ``empirically
complete'' for the projective sets and $ZFC+PD$ ``rigidifies'' the
properties of projective sets w.r.t. forcing: it makes them
``forcing-absolute'' or ``generically absolute''. One can also consider the
even stronger axiom (Woodin axiom) ``$L(\mathbb{R)}$ satisfies $AD$'' where $L(%
\mathbb{R)}$ is the constructible closure of $\mathbb{R}$ (i.e. the smallest inner
model containing $On$ and $\mathbb{R}$, see above) and the \emph{Axiom of
Determinacy} $AD$ means that \emph{every} subset of $\mathbb{R}$ is determined. 
$AD$ is incompatible with $AC$ since, according to Fubini theorem, $AC$
enables to construct non-Lebesgue measurable, and therefore non-determined,
subsets of $\mathbb{R}$.

\section{The transcendance of $\mathbb{R}$ over $L$ and the set $0^{\#}$ ($0$
sharp)}

\subsection{Indiscernible ordinals}

Once we accept the relevance and the legitimacy of $LCAs$, we need some
tools for measuring the transcendence of $V$ over $L$. A first possibility
is given by what are called \emph{indiscernible} ordinals (Silver, $1966$)
which enable to construct the simplest canonical non constructible real. We
consider the levels of $L$ of the form $\left\langle L_{\lambda },\in
\right\rangle $ with $\lambda $ a limit ordinal. A set $I$ of ordinals in
this cumulative hierarchy $L_{\lambda }$ of constructive sets up to level $%
\lambda $ is called a set of indiscernibles if, for every $n$-ary formula $%
\varphi \left( x_{1},\ldots ,x_{n}\right) $, the validity of $\varphi $ on $%
I $ is independent of the choice of the $x_{i}$'s: that is for every
sequences $c_{1}<\ldots <c_{n}$ and $d_{1}<\ldots <d_{n}$ in $I$%
\[
L_{\lambda }\models \varphi \left( c_{1},\ldots ,c_{n}\right) \text{ iff }
L_{\lambda }\models \varphi \left( d_{1},\ldots ,d_{n}\right) . 
\]

When it exists, the set $S$ of \emph{Silver indiscernibles} is characterized
by the following properties, which express that, for all \emph{uncountable}
cardinals $\kappa $, all the $L_{\kappa }$'s share essentially the same
first-order structure:

\begin{enumerate}
\item  $\kappa \in S$ (all uncountable cardinals of $V$ are indiscernible in 
$L$).

\item  $S\cap \kappa $ is of order-type $\kappa $.

\item  $S\cap \kappa $ is closed and unbounded (club, see above) in $\kappa $
if $\kappa $ is regular.

\item  $S\cap \kappa $ is a set of indiscernibles for $\left\langle
L_{\kappa },\in \right\rangle $.

\item  The \emph{Skolem hull} of $S\cap \kappa $ in $L_{\kappa }$ is equal
to $L_{\kappa }$: $\func{Hull}^{L_{\kappa }}\left( S\cap \kappa \right)
=L_{\kappa }$, where the Skolem hull of $I\subset L_{\kappa }$ is
constructed by adding for every $(n+1)$-ary formula $\varphi \left(
y,x_{1},\ldots ,x_{n}\right) $ with $x_{i}\in I$ a Skolem term $t_{\varphi
}\left( x_{1},\ldots ,x_{n}\right) $ which is the smallest $y$ (for the
wellorder of $L$) s.t. $\varphi \left( y,x_{1},\ldots ,x_{n}\right) $ if
such an $y$ exists and $0$ otherwise. In other words every constructible
element $a\in L_{\kappa }$ is \emph{definable} by a definite description
with parameters in the indiscernibles $S\cap \kappa $.
\end{enumerate}

This can be generalized to structures $\mathcal{M}=\left\langle
M,E\right\rangle $ with a binary relation $E$ looking like $\in $, that is,
which are elementary equivalent to some $\left\langle L_{\lambda },\in
\right\rangle $ for $\lambda $ a limit ordinal and $I\subset M$. In that
case, we have $\func{Hull}^{\mathcal{M}}\left( I\right) \prec \mathcal{M}$
and in fact $\func{Hull}^{\mathcal{M}}\left( I\right) $ is the smallest
elementary substructure of $\mathcal{M}$ containing $I$. Let $\Sigma =\Sigma
\left( \mathcal{M},I\right) $ be the set of formulae $\varphi $ which can be
satisfied by $\mathcal{M}$ on $I$.\ This defines particular sets of formulae
called $EM$-sets (from Ehrenfeucht-Mostowski, $1956$). The EM theorem says
that if $\Sigma $ is a theory having infinite models and if $\left\langle
I,<\right\rangle $ is any total well-ordering of infinite order-type $\alpha
\geq \omega $, then there exists a model $\mathcal{M}$ of $\Sigma $
containing $I$ for which $I$ is a set of indiscernibles, and moreover, $%
\mathcal{M}$ can be chosen in such a way as to be the \emph{Skolem hull} of $%
I$: $\mathcal{M}=\func{Hull}^{\mathcal{M}}\left( I\right) $. Such an $\left( 
\mathcal{M},I\right) $ is essentially unique and its transitive collapse
(isomorphism with a structure where $E$ becomes $\in $) is written $\left( 
\mathcal{M}\left( \Sigma ,\alpha \right) ,I\left( \Sigma ,\alpha \right)
\right) $ where $\alpha $ is the order type of $I$ ($I\left( \Sigma ,\alpha
\right) $ is therefore a set of true $\in $-ordinals).

One can develop a theory of $EM$-sets and of their well-foundedness. If $%
\Sigma $ is well-founded (i.e. if $\mathcal{M}\left( \Sigma ,\alpha \right) $
is well-founded for every ordinal $\alpha $) and if $\alpha $ is a limit
ordinal, then $\mathcal{M}\left( \Sigma ,\alpha \right) $ is isomorphic to a 
$\left\langle L_{\lambda },\in \right\rangle $.\ Moreover, if $I\left(
\Sigma ,\kappa \right) $, with $\kappa >\omega $ an uncountable cardinal, is
unbounded in the class of ordinals of $\mathcal{M}\left( \Sigma ,\alpha
\right) $ (and it is then the case for every ordinal $\alpha >\omega $), and
if for every ordinal $\gamma <i_{\omega }$ (the $\omega $-th element of $%
I\left( \Sigma ,\kappa \right) $) we have $\gamma \in \func{Hull}^{\mathcal{M%
}}\left( \left\{ i_{n}\right\} \right) $ (and it is then the case for every
ordinal $\alpha >\omega $), then $\mathcal{M}\left( \Sigma ,\kappa \right)
=\left\langle L_{\kappa },\in \right\rangle $, $I\left( \Sigma ,\kappa
\right) $ is closed unbounded in $\kappa $ and if $\tau >\kappa >\omega $
then $I\left( \Sigma ,\tau \right) \cap \kappa =I\left( \Sigma ,\kappa
\right) .$

Let us return to Silver indiscernibles. If such a $\Sigma $ exists, $S$ is
defined by 
\[
S=\bigcup \left\{ I\left( \Sigma ,\kappa \right) :\kappa \text{ uncountable
cardinal}\right\} . 
\]

The uniqueness of $S$ is a consequence of the unicity of such a $\Sigma $:

\textbf{Theorem.}\ \textsl{Such a }$\Sigma $\textsl{\ is unique and is the
set of }$n$\textsl{-ary sentences }$\varphi $\textsl{\ s.t. }$L_{\aleph
_{\omega }}\models \varphi \left( \aleph _{1},\ldots ,\aleph _{n}\right) $%
\textsl{.\ It is called ``zero sharp'' and written }$0^{\#}$ \textsl{(see
below).\hfill }$\square $

If there exists an uncountable limit cardinal $\kappa $ s.t. $\left\langle
L_{\kappa },\in \right\rangle $ possesses an uncountable set $I$ of
indiscernibles, then $S$ exists. The existence of $S$ is also implied by
large cardinal hypothesis as for instance:

\textbf{Theorem}.\ \textsl{If there exists a }$MC$\textsl{\ then }$S$\textsl{%
\ exists and moreover }$L_{\kappa }\prec L_{\lambda }$\textsl{\ for every
uncountable }$\kappa <\lambda $\textsl{.\hfill }$\square $

The existence of $S$ under $MC$ means that after the first uncountable level 
$L_{\aleph _{1}}$ all the $L_{\kappa }$ share essentially the same
first-order theory.\ $V$ transcends $L$, but in such a way that it makes $L$
as simple as possible, $L_{\aleph _{1}}=\func{Hull}^{L_{\aleph _{1}}}\left(
S\cap \aleph _{1}\right) \ $determining the theory of $L$.

A deep consequence is that the \emph{truth} in $L$ becomes definable in $V$.
Indeed, let $\varphi \left( x_{1},\ldots ,x_{n}\right) $ be a formula. There
exists an uncountable cardinal $\kappa $ s.t. 
\[
\text{for all }\left( x_{i}\right) \in L_{\kappa }\text{, }L\models \varphi
\left( x_{i}\right) \text{ iff }L_{\kappa }\models \varphi \left(
x_{i}\right) \text{ .} 
\]

\noindent As $L_{\kappa }\prec L_{\lambda }$ if $\kappa <\lambda $, we have 
\[
L\models \varphi \left( x_{i}\right) \text{ iff }L_{\lambda }\models \varphi
\left( x_{i}\right) \text{ for all }\lambda \geq \kappa \;. 
\]

\noindent Now, we arithmetize the situation.\ Let $T=\left\{ \ulcorner
\varphi \urcorner :L_{\aleph _{1}}\models \varphi \right\} $ be the set of
G\"{o}del numbers of the $\varphi $ valid in $L_{\aleph _{1}}$ and therefore
in all the $L_{\kappa }$ ($\kappa $ uncountable) by elementarity.\ Then 
\[
L\models \varphi \text{ iff }\ulcorner \varphi \urcorner \in T 
\]

\noindent defines the truth in $L$.\ This is not in contradiction with
G\"{o}del-Tarski uncompleteness theorems since $\aleph _{1}$ and $T$ are 
\emph{not definable} in $L$ and therefore the truth of $L$ is not definable
in $L$.

\subsection{The set $0^{\#}$}

As $L_{\aleph _{\omega }}\prec L$ and $\aleph _{i}\in S$ for $i>0$ in $%
\omega $, we can represent the indiscernibles in formulae by some of the $%
\aleph _{i}$'s and restrict to $L_{\aleph _{\omega }}$, which contains all
the $\aleph _{i}$'s. Then, $L\models \varphi \left( x_{i}\right) $ for $%
x_{i}\in S$ iff $L_{\aleph _{\omega }}\models \varphi \left( \aleph
_{i}\right) $. Solovay called $0^{\#}$ (zero sharp) the set (if it exists)
defined by 
\[
0^{\#}=\left\{ \varphi :L_{\aleph _{\omega }}\models \varphi \left( \aleph
_{i}\right) \right\} 
\]

\noindent which is the set of formulae true on the indiscernibles of $L$.\
Via G\"{o}delization $0^{\#}$ becomes a set of integers (also written $%
0^{\#} $) and can therefore be coded by a real (also written $0^{\#}$).

We must emphasize the fact that, as $L_{\aleph _{1}}\prec L$, every
constructible set $x\in L$ which is \emph{definable} in $L$ is \emph{%
countable} since its definite description is valid in $L_{\aleph _{1}}$ by
elementarity and therefore $x\in L_{\aleph _{1}}$. More generally for every
infinite constructible set $x\in L$ we have $\left| \mathcal{P}(x)\right|
^{L}=\left| x\right| $. Since the existence of a measurable cardinal implies
that $0^{\#}$ exists, we have:

\textbf{Corollary}.\ \textsl{If there exists a }$MC$\textsl{, the
constructible continuum }$\mathbb{R}^{L}$\textsl{\ is countable.\hfill }$%
\square $

Via arithmetization through G\"{o}del numbers, the non constructible set $%
0^{\#}$ can be considered as a very special subset $\notin L$ of $\omega =%
\mathbb{N}$ or as a very special real number coding the truth in $L.$ Its
existence implies that every uncountable cardinal $\kappa $ of $V$ is an
indiscernible of $L$ and shares \emph{all} large cardinal axioms verified by 
$L$.

A property equivalent to the existence of $0^{\#}$ is the \emph{non rigidity}
of $L$:

\textbf{Theorem (Kunen).}\ $0^{\#}$ exists iff there exists a non trivial
elementary embedding $j:L\prec L$ (this presuppose $V\neq L$ and $j$ non
trivial, see below).$\hfill \square $

Indeed, as $\func{Hull}^{L}(S)=L$, for every $x\in L$ there exists a Skolem
term $t$ s.t. $x=t\left( i_{\alpha _{1}},\ldots ,i_{\alpha _{n}}\right) $, $%
i_{\alpha }$ being the $\alpha $-th element of $S$. $j$ is then simply
defined by \emph{the shift} on indiscernibles 
\[
j(x)=j\left( t\left( i_{\alpha _{1}},\ldots ,i_{\alpha _{n}}\right) \right)
=t\left( i_{\alpha _{1}+1},\ldots ,i_{\alpha _{n}+1}\right) . 
\]

\noindent One shows that it is an elementary embedding and, as $j\left(
i_{0}\right) =i_{1}\neq i_{0}$, $j$ is non trivial. $\hfill \square $

By the way, this proves again that $V\neq L$ since another theorem of Kunen
proves that 
\[
ZFC\vdash \text{ there exists no }j:V\prec V\;. 
\]

The existence of $0^{\#}$ is a principle of transcendence of $V$ over $L\;$%
expressing that $V$ is very different from $L$.\ If $0^{\#}$ doesn't exist,
then $V$ looks like $L$ ($L$ is a good approximation of $V$) according to
the result:

\textbf{Covering lemma (Jansen)}.\ \textsl{If }$0^{\#}$\textsl{\ doesn't
exist, then if }$x$\textsl{\ is an uncountable set of ordinals there exists
a constructible set }$y\supseteq x$\textsl{\ of the same cardinality as }$x$%
\textsl{. So, every set }$x$\textsl{\ of ordinals can be covered by a
constructible set }$y\supseteq x$\textsl{\ of cardinality }$\left| y\right|
=\left| x\right| .\aleph _{1}$.\textsl{\hfill }$\square $

\textbf{Corollary.} \textsl{If }$0^{\#}$\textsl{\ doesn't exist, the
covering lemma implies that, for every limit singular cardinal }$\kappa $ 
\textsl{of }$V$\textsl{, we have }$\left( \kappa ^{+}\right) ^{L}=\kappa
^{+} $\textsl{, which shows that }$V$\textsl{\ and }$L$\textsl{\ are quite
similar.} \textsl{\hfill }$\square $

Indeed (see Jech [1978], p.\ 358), if $\lambda =\left( \kappa ^{+}\right)
^{L}$ and if $\lambda <\kappa ^{+}$, then $\left| \lambda \right| =\kappa $
.\ But as $\kappa $ is singular, we would have $\func{cf}\left( \lambda
\right) <\left| \lambda \right| $ and this is impossible since $\lambda $ is
regular in $L$ and $\lambda \geq \omega _{2}$. For, if $x$ is an unbounded
subset of $\lambda $ of cardinal $\left| x\right| =\func{cf}\left( \lambda
\right) ,$ it can be covered by a \emph{constructible} subset $y\in L$ of $%
\lambda $ of cardinal $\left| y\right| =\left| x\right| .\aleph _{1}$ and,
as $\lambda $ is regular in $L$, $\left| y\right| =\left| \lambda \right| $.
So, $\left| \lambda \right| =\aleph _{1}.\func{cf}\left( \lambda \right) $
and, as $\lambda \geq \omega _{2}$, $\left| \lambda \right| =\func{cf}\left(
\lambda \right) $.$\hfill \square $

\textbf{Corollary.} \textsl{If }$GCH$\textsl{\ fails at a strong limit
singular cardinal, then }$0^{\#}$\textsl{\ exists.}$\hfill \square $

Indeed, if $0^{\#}$ doesn't exist, for such a singular cardinal $\kappa $ we
have $\left( \kappa ^{+}\right) ^{L}=\kappa ^{+}$.\ As $L$ satisfies $GCH$ , 
$\left( 2^{\kappa }\right) ^{L}=\left( \kappa ^{+}\right) ^{L}=\kappa ^{+}$.%
\footnote{%
As, if $0^{\#}$ exists, $\left( 2^{\kappa }\right) ^{L}=\kappa $ for every
infinite cardinal since $\kappa ^{+}$ is inaccessible in $L$, $\left(
2^{\kappa }\right) ^{L}=\kappa ^{+}>\kappa $ is a counter exemple when $%
0^{\#}$ doesn't exist.} Now, $\kappa $ is a strong limit by hypothesis (i.e. 
$\lambda <\kappa \Rightarrow 2^{\lambda }<\kappa $) and this implies $\kappa
^{\func{cf}\left( \kappa \right) }=2^{\kappa }$ and moreover, since $\kappa $
is singular and therefore $\func{cf}\left( \kappa \right) <\left| \kappa
\right| $, $2^{\func{cf}\left( \kappa \right) }<\kappa $. Let $x\in A=\left[
\kappa \right] ^{\func{cf}\left( \kappa \right) }$ be a subset of $\kappa $
of cardinal $\func{cf}\left( \kappa \right) $.\ It is covered by a
constructible subset $y\in L$ of $\kappa $ of cardinal $\left| y\right|
=\lambda =\aleph _{1}.\func{cf}\left( \kappa \right) $.\ Then $A$ can be
covered by the union of the $\left[ y\right] ^{\func{cf}\left( \kappa
\right) }$ for such $Y$.\ But $\left| \left[ y\right] ^{\func{cf}\left(
\kappa \right) }\right| =\lambda ^{\func{cf}\left( \kappa \right) }=\left(
\aleph _{1}.\func{cf}\left( \kappa \right) \right) ^{\func{cf}\left( \kappa
\right) }=2^{\func{cf}\left( \kappa \right) }$, and by hypothesis $2^{\func{%
cf}\left( \kappa \right) }<\kappa $.\ Now, there exist at most $\left| 
\mathcal{D}\left( \kappa \right) \right| $ such $y$~\footnote{%
Recall that $\mathcal{D}(x)$ is the set of constructible parts of $x$ (see
above).} and $\left| \mathcal{D}\left( \kappa \right) \right| =\left( \kappa
^{+}\right) ^{L}=\kappa ^{+}$.\ All this implies $\left| A\right| =\kappa ^{%
\func{cf}\left( \kappa \right) }=\kappa ^{+}$.\ As $\kappa ^{\func{cf}\left(
\kappa \right) }=2^{\kappa }$, we have $2^{\kappa }=\kappa ^{+}$, that is $%
GCH$. $\hfill \square $

There are many ways for insuring that $0^{\#}$ exists, for instance $PFA$
since $PFA$ implies the failure of $\square _{\kappa }$ for every $\kappa $
(Todorcevic, 1980). More generally, many results show that a failure of $%
\square _{\kappa }$ is linked with $LCAs$.\ For instance:

\begin{itemize}
\item  Jensen: if $\square _{\kappa }$ fails for some singular $\kappa $,
there exists an inner model $M$ with a strong cardinal.

\item  Solovay: if $\kappa $ is supercompact, then $\square _{\lambda }$
fails for every $\lambda >\kappa $.
\end{itemize}

\noindent This is related to the fact that properties such as the covering
between $V$ and $M$ allow to reflect $\square _{\kappa }$ from $M$ to $V$
and therefore, if $\square _{\kappa }$ fails in $V$ for $\kappa $ singular
it is because there exist some $LCs$ violating the covering lemma.

When $0^{\#}$ exists, a very interesting structure to look at is $L\left[
0^{\#}\right] $.\ It can be shown~\footnote{%
See e.g. Steel [2001].} that $L\left[ 0^{\#}\right] $ has a fine structure,
satisfies the Global Square property $\forall \kappa \;\square _{\kappa }$
and, for every singular cardinal $\kappa $ of $V$, $\left( \kappa
^{+}\right) ^{L\left[ 0^{\#}\right] }=\kappa ^{+}$.\ Iterating the sharp
operation, one can get an increasing sequence of models $M$ mildly
transcendent over $L$ which have also a fine structure, satisfy the Global
Square property and, for every singular cardinal $\kappa $ of $V$, $\left(
\kappa ^{+}\right) ^{M}=\kappa ^{+}$, and are smaller than the first inner
model possessing a measurable cardinal. Under the hypothesis that $a^{\#}$
exists for every set of ordinals $a$ ($a^{\#}$ is defined as $0^{\#}$ but in 
$L[a]$), one can even go beyond this limit, up to the existence of a
supercompact cardinal.

\subsection{$0^{\#}$ and the hierarchical structure of $V$ beyond $L$.}

The equivalence between the existence of $0^{\#}$ and the existence of a non
trivial elementary embedding $j:L\prec L$ enables to clarify the structure
of $V$ beyond $L$~\footnote{%
See e.g. Schimmerling [2001].}. Let $\kappa =\func{crit}\left( j\right) $
and let $\mathcal{U}$ be the set of subsets $u\subseteq \kappa $ s.t. $u\in
L $ (i.e. $u\in L\cap \mathcal{P}\left( \kappa \right) $) and $\kappa \in
j\left( u\right) $. $\mathcal{U}$ is trivially a filter.\ It is an
ultrafilter since for every $u\in L\cap \mathcal{P}\left( \kappa \right) $
either $\kappa \in j\left( u\right) $ or $\kappa \notin j\left( u\right) $.
It is a free ultrafilter since if $u\subset \kappa $ is bounded, then $%
\kappa \notin u$, $j\left( u\right) =u$, $\kappa \notin j\left( u\right) $
and therefore $u\notin \mathcal{U}$. Moreover, it is $\kappa $-complete
w.r.t. $L$ since if $u_{\alpha }\in \mathcal{U}\cap L$ is a family with $%
\alpha <\beta <\kappa $ then $\bigcap_{\alpha <\beta }u_{\alpha }\in 
\mathcal{U}$.\ One can show that the ultrapower $L^{\mathcal{U}}$ is
well-founded. Due to \L o\v {s} theorem, the embedding $i:L\prec L^{\mathcal{%
U}}$ is elementary.\ But the Mostowski collapsing lemma implies that $%
\left\langle L^{\mathcal{U}},\in _{\mathcal{U}}\right\rangle \simeq
\left\langle M_{\mathcal{U}},\in \right\rangle $ for some transitive inner
model $M_{\mathcal{U}}$.\ But necessarily $M_{\mathcal{U}}=L$ by minimality
of $L$ and via this isomorphism $i:L\prec L^{\mathcal{U}}$ becomes an
elementary embedding $j:L\prec L$. It can be shown that, if $\lambda =\left(
\kappa ^{+}\right) ^{L}$ then $\mathcal{M}=\left\langle L_{\lambda },\in ,%
\mathcal{U}\right\rangle $ is a model of $ZF-$\{Power Set axiom\} where $%
\kappa $ becomes the largest cardinal, $\mathcal{U}$ remains a free
ultrafilter $\kappa $-complete with ``good'' technical properties
(``normality'' and ``amenability'').\ Such a procedure can be \emph{iterated 
}on the ordinals. Starting from a $\mathcal{M}_{0}=\left\langle L_{\lambda
_{0}},\in ,\mathcal{U}_{0}\right\rangle $, one gets a $\mathcal{M}%
_{1}=\left\langle L_{\lambda _{1}},\in ,\mathcal{U}_{1}\right\rangle $, etc.
The successive $\mathcal{M}_{\alpha }$ yield a sequence of critical
cardinals $\kappa _{\alpha }$ which are indiscernibles for $L$.

\section{Determination and reflection phenomena}

To measure the size of large cardinals, the best way is to use associated 
\emph{reflection }phenomena which are of a very deep philosophcal value.%
\footnote{%
See Martin-Steel [1989], Patrick Dehornoy [1989].} Intuitively, reflection
means that the properties of the whole universe $V$ are reflected in
sub-universes. As was emphasized by Matthew Foreman (1998, p.\ 6):

\begin{quotation}
\noindent ``Any property that holds in the mathematical universe should hold
of many set-approximations of the mathematical universe.''
\end{quotation}

\textbf{Definition.} \textsl{A cardinal} $\chi $ reflects \textsl{a relation 
}$\Phi (x,y)$\textsl{\ defined on ordinals if every solution }$y\geq \chi $ 
\textsl{parametrized by }$x<\chi $\textsl{\ can be substituted for by a
solution }$y<\chi $:

\begin{center}
$\forall \alpha (\in On)<\chi \ \left[ \exists \beta \geq \chi \ \Phi
(\alpha ,\beta )\Rightarrow \exists \beta ^{*}<\chi \ \Phi (\alpha ,\beta
^{*})\right] .$
\end{center}

Let $j$ be an elementary embedding $j:M\prec M^{*}$. $\chi =\func{crit}%
\left( j\right) $ is a large~-- in fact at least measurable~-- cardinal,
which increases indefinitely when $M^{*}$ moves near to $M$, the limit $%
M^{*}=M$ being inconsistent according to Kunen theorem.

To see that it is a reflection phenomenon, let $\Phi (\alpha ,\chi )$ be a
relation that holds in $M$ for $\alpha <\chi $. If $M^{*}$ is sufficiently
close to $M$ for $\Phi (\alpha ,\chi )$ to remain true in $M^{*}$, then $%
M^{*}\models \exists (x<j(\chi ))\Phi (\alpha ,x)$ (it is sufficient to take 
$x=\chi $). But, according to the elementarity of the embedding $j$, this is
equivalent to $M\models \exists (x<\chi )\Phi (\alpha ,x)$.

To go beyond measurable cardinals, specialists use the following technique.
Let $V_{\alpha }$ be the cumulative hierarchy of sets up to level $\alpha $.
For $\chi $ critical (and therefore measurable), one has $V_{\chi
}^{M^{*}}=V_{\chi }^{M}$ (that is the equality of $M$ and $M^{*}$ up to
level $\chi $).

\textbf{Definition.} \textsl{The cardinal }$\chi $\textsl{\ is called} \emph{%
\ \ superstrong} \textsl{in }$M$\textsl{\ if there exists an elementary
embedding }$j$\textsl{\ s.t. }$V_{j(\chi )}^{M^{*}}=V_{j(\chi )}^{M}$ 
\textsl{(that is }$V_{j(\chi )}^{M^{*}}\subset M$\textsl{\ and }$M=M^{*}$%
\textsl{\ up to }$j(\chi )$\textsl{\ and not only up to }$\chi $\textsl{\
).\hfill }$\square $

Between measurable and superstrong cardinals, Hugh Woodin introduced another
class of large cardinals.

\textbf{Definition.} \textsl{A cardinal }$\delta $\textsl{\ is called a
Woodin cardinal if for every map }$F:\delta \rightarrow \delta $\textsl{,
there exists} $\kappa <\delta $\textsl{\ and an elementary embedding }$j$%
\textsl{\ of critical ordinal }$\kappa $\textsl{\ s.t. }$F|_{\kappa }:\kappa
\rightarrow \kappa $\textsl{\ and }$V_{j(F(\kappa ))}^{M^{*}}=V_{j(F(\kappa
))}^{M}$\textsl{\ (that is }$M=M^{*}$\textsl{\ up to }$j(F(\kappa ))$\textsl{%
\ ).\hfill }$\square $

Woodin has shown that:

\begin{enumerate}
\item  if $\delta $ is a Woodin cardinal, there exist infinitely many
smaller measurable cardinals $\chi <\delta $,

\item  if $\lambda $ is a superstrong cardinal, there exist infinetely many
smaller Woodin cardinals $\delta <\lambda $.
\end{enumerate}

A key result is the Martin-Steel theorem which evaluates exactly the
``cost'' of determinacy:

\textbf{Martin-Steel theorem (1985).} \textsl{If there exist }$n$ \textsl{%
Woodin\ cardinals }$\delta _{i}$\textsl{,\ }$i=1,\ldots ,n$\textsl{,} 
\textsl{dominated by a measurable cardinal }$\kappa $ \textsl{(}$\kappa
>\delta _{i}$\textsl{\ for all} $i$\textsl{), then }$ZFC\vdash \func{Det}
\left( \Pi _{n+1}^{1}\right) .\hfill \square $

The converse is due to Woodin.

\textbf{Corollary.} \textsl{If there}$\ $\textsl{exists a countable infinity
of Woodin cardinals dominated by a measurable cardinal, in particular if
there exists a superstrong cardinal }$\lambda $,\textsl{\ then }Projective
Determinacy\textsl{\ is valid (all the projective subsets of }$\mathbb{R}$ 
\textsl{are determined).\hfill }$\square $

Projective Determinacy is also valid under $PFA$ (Woodin, see above).

It is for this reason that specialists consider that $ZFC+$ Projective
Determinacy is a ``good'' axiomatic for $\mathbb{R}$. We must also emphasize
the:

\textbf{Martin-Steel-Woodin theorem (1987).}\textsl{\ If there}$\ $\textsl{%
exists a countable infinity of Woodin cardinals dominated by a measurable
cardinal, in particular if there}$\ $\textsl{exists a superstrong cardinal }$%
\lambda $,\ \textsl{then }$L(\mathbb{R})$\textsl{\ (the smallest inner model of 
}$V$\textsl{\ containing the ordinals }$On$\textsl{\ and }$\mathbb{R}$\textsl{,
see above) satisfies the axiom of complete determinacy AD: every} $%
A\subseteq \mathbb{R}$\textsl{\ is determined. (This result is stronger than
the previous one since }$\mathcal{P}(\mathbb{R})\cap L(\mathbb{R})$ \textsl{is a
larger class than the projective class.)\hfill }$\square $

$AD$ is incompatible with $AC$ since $AC$ enables the construction of a non
determined well ordering on $\mathbb{R}$ (see above).\footnote{%
So, the inner model $L(\mathbb{R})$ of a $ZFC$--model $V$ can violate $AC$.}

But the most significative results concern perhaps the situation where no
property of $\mathbb{R}$ could be further modified in a forcing extension. In
that case, the theory of the continuum becomes ``\emph{rigid}''. Woodin and
Shelah have shown that it is possible to approximate this ideal goal if
there exists a \emph{supercompact }cardinal $\kappa $. $\kappa $ is $\gamma $%
-supercompact if there exists an elementary embedding $j:V\prec M$ s.t. $%
\func{crit}(j)=\kappa $, $\gamma <j(\kappa )$ and $M^{\gamma }\subseteq M$.\ 
$\kappa $ is supercompact if it is $\gamma $-supercompact for every $\gamma
\geq \kappa $ ($\kappa $ is $\kappa $-supercompact iff it is measurable).%
\footnote{%
See Dehornoy [2003].}\ Such a deep result clarifies the nature of the axioms
which are needed for a ``good'' theory of the continuum.

\section{Woodin's $\Omega $-logic}

Large cardinal axioms ($LCAs$) can decide some properties of regularity of $%
\mathbb{R}$, but they cannot settle $CH$ since a ``small'' forcing (adding $%
\aleph _{2}$ new subsets to $\omega $) is sufficient to force $\lnot CH$
from a $CH$-model and such a small forcing remains possible irrespective of
what $LCAs$ are introduced (Levy-Solovay theorem). We need therefore a new
strategy.\ As we have just seen, the most natural one is to try to make the
properties of the continuum \emph{\ immune} relatively to forcing, that is
to make the continuum in some sense ``\emph{rigid}''. The deepest
contemporary results in this perspective are provided by Woodin's
extraordinary recent works on $\Omega $-logic and the negation of $CH$.

We look for theories sharing \emph{some absoluteness properties relatively
to forcing}.\ This is called ``\emph{conditional generic absoluteness}''.%
\footnote{%
See Steel [2004]: ``Generic absoluteness and the continuum problem''.}

The fragment of $V$ where $CH$ ``lives'' naturally is $\left( H_{2},\in
\right) $ where $\left( H_{k},\in \right) $ is the set of sets $x$ which are
hereditary of cardinal $\left| x\right| <\aleph _{k}$. The fragment $\left(
H_{0},\in \right) =V_{\omega }$ is the set of hereditary finite sets and,
with the axioms $ZF$ minus the axiom of infinity, is equivalent to first
order arithmetic $\left\langle \omega =\mathbb{N},+,.,\in \right\rangle $ with
Peano axioms.\ In one direction, $\mathbb{N}$ can be retrieved from $H_{0}$
using von Neumann's construction of ordinals and, conversely, $H_{0}$ can be
retrieved from $\mathbb{N\,}$via Ackermann's trick: if $p,q$ are integers, $%
p\in q$ iff the $p$-th digit in the binary extension of $q$ is $1$. For
first order arithmetic, Peano axioms are ``empirically'' and practically
complete in spite of G\"{o}del incompleteness theorem. The following
classical result expresses their ``rigidity'':

\textbf{Sch\"{o}nfield theorem}.\ $H_{0}$\textsl{\ is absolute and a
fortiori forcing-invariant. Incompleteness cannot be manifested in it using
forcing.\hfill }$\square $

We can therefore consider $ZFC$ as a ``good'' theory for first order
arithmetic. But it is no longer the case for larger fragments of $V$.

The fragment $\left( H_{1},\in \right) $ of $V$ composed of countable sets
of finite ordinals is isomorphic to $\left\langle \mathcal{P}\left( \omega
\right) =\mathbb{R},\omega ,+,.,\in \right\rangle $ and corresponds to second
order arithmetic (i.e. analysis).\ The definable subsets $A\subseteq 
\mathcal{P}\left( \omega \right) $ are the projective subsets~\footnote{$%
A\subseteq \mathcal{P}\left( \omega \right) $ is definable in $H_{1}$ (with
parameters in $H_{1}$) iff there exists a first-order formula $\varphi (x,y)$
and a parameter $b\in \mathcal{P}\left( \omega \right) $ s.t. $A=\left\{
a\in \mathcal{P}\left( \omega \right) \left| H_{1}\models \right. ~\varphi
(a,b)\right\} $.\ If $\pi :\mathcal{P}\left( \omega \right) \rightarrow
[0,1]\simeq \mathbb{R}$ is given by $\pi (a)=\sum_{i\in a}\frac{1}{2^{i}}$,
then $X\subseteq [0,1]$ is projective iff $A=\pi ^{-1}(X)$ is definable in $%
H_{1}$.} and therefore $H_{1}$ can be considered as the fragment of $V$
where the projective sets live. We have seen that to settle and ``freeze''
most of its higher order properties (regularity of projective sets) w.r.t.
forcing, we need $LCAs$ and in particular $PD$.

As is emphasized by Woodin (2003, quoted in Dehornoy [2003]):

\begin{quotation}
\noindent ``\emph{Projective Determinacy} settles (in the context of $ZFC$)
the classical questions concerning the projective sets and moreover Cohen's
method of forcing \emph{cannot} be used to establish that questions of
second order number theory are formally unsolvable from this axiom. (...) I
believe the axiom of \emph{Projective Determinacy} is as true as the axioms
of Number Theory.\ So I suppose that I advocate a position that might best
be described as \emph{Conditional Platonism}.''
\end{quotation}

We have also seen that under the $LCA$ ``there exists a proper class of
Woodin cardinals'' ($PCW$: for every cardinal $\kappa $ there exists a
Woodin cardinal $>\kappa $) we have:

\textbf{Theorem (Woodin, 1984)}.\ $ZFC+PCW\vdash H_{1}$ \textsl{is immune
relatively to forcing in the sense its properties are
forcing-invariant.\hfill }$\square $

As $PCW$ implies at the same time $PD$ and forcing-invariance for $H_{1}$,
it can be considered as a ``good'' theory, ``empirically'' and
``practically'' complete (marginalizing incompleteness) for $\left(
H_{1},\in \right) $, that is for analysis (second order arithmetic). $PCW$
implies the \emph{generic completeness} result that all the $L[\mathbb{R}]$ of
generic extensions $V[G]$ are elementary equivalent.

The idea is then to try to generalize these types of absoluteness properties
relative to forcing.\ The general strategy for deciding that way $ZFC$%
-undecidable properties $\varphi $ in a fragment $H$ of $V$ is described in
the following way by Patrick\ Dehornoy (2003):

\begin{quotation}
\noindent ``every axiomatization freezing the properties of $H$ relatively
to forcing (i.e. neutralizing forcing at the level $H$) implies $\varphi $''.
\end{quotation}

The main problem tackled by Woodin was to apply this strategy to the
fragment $\left( H_{2},\in \right) $\ of $V$ which is associated to the set $%
\mathcal{P}\left( \omega _{1}\right) $ of countable ordinals. $\mathcal{P}%
\left( \omega _{1}\right) $ is not $\mathcal{P}\left( \mathbb{R}\right) $ if $%
\lnot CH$ is satisfied, but nevertheless it is possible to code $CH$ by an $%
H_{2}$-formula $\varphi _{CH}$ s.t. $H_{2}\models \varphi _{CH}$ is
equivalent to $CH$.\footnote{%
The point is rather technical. Woodin has shown that if $\lnot CH$ is valid
(i.e. $\mathbb{R}>\omega _{1}$), $\mathcal{P}(\mathbb{R})\notin H_{2}$ and is
already too big for freezing (neutralizing the effects of forcing) the
fragments of $V$ containing it.} The problem with $H_{2}$ is that ``small''
forcings preserve $LCAs$~\footnote{%
Large cardinal axioms are axioms of the form $A=\exists \kappa \psi (\kappa
) $ which share the property that if $V\vDash A$ then the cardinal $\kappa $
is inaccessible and $\psi (\kappa )$ is forcing-invariant for every forcing
extension $V[G]$ of forcing cardinal $<\kappa $ (``small'' forcings).} and
in particular (Levy-Solovay theorem, 1967) a small forcing of cardinal $%
\aleph _{2}$ that enables to violate $CH$ by adding $\aleph _{2}$ subsets to 
$\mathbb{N}$ preserves $LCAs$. Therefore $H_{2}$ \emph{cannot} be rigidified by 
$LCAs$.\footnote{%
See Dehornoy (2003).} Whatever the large cardinal hypothesis $A$ may be,
there will be always generic extensions $M$ and $N$ of $V$ both satisfying $%
A $ such that $M\models CH$ and $N\models \lnot CH$.\ As $CH$ is equivalent
to a $\Sigma _{2}^{1}$ formula, $M$ and $N$ cannot be elementary equivalent
from the $\Sigma _{2}^{1}$ level.

Woodin's fundamental idea to overcome the dramatic difficulties of the
problem at the $H_{2}$ level was to \emph{strengthen logic} by restricting
the admissible models and constructing a new logic adapted to
forcing-invariance or ``\emph{generic invariance}''. As he explains in his
key paper on ``The continuum hypothesis'' (2001, p.~682):

\begin{quotation}
\noindent ``As a consequence (of generic invariance), any axioms we find
will yield theories for $\left\langle H\left[ \omega _{2}\right] ,\in
\right\rangle ,$ whose `completeness' is immune to attack by applications of
Cohen's method of forcing, just as it is the case for number theory.''
\end{quotation}

In a first step, he introduced the notion of $\Omega $-validity $\models
_{\Omega }$ also called in a first time $\Omega ^{*}$-derivability $\vdash
_{\Omega ^{*}}$.

\textbf{Definition}.\ $T$\textsl{\ being a theory in }$ZFC$\textsl{, we have 
}$T\models _{\Omega }\varphi $\textsl{\ iff }$\varphi $\textsl{\ is valid in
every generic extension where }$T$\textsl{\ is valid, that is iff for every
generic extension }$V\left[ G\right] $\textsl{\ and every level }$\alpha $%
\textsl{, }$\left( V_{\alpha }\right) ^{V\left[ G\right] }\models T$\textsl{%
\ implies }$\left( V_{\alpha }\right) ^{V\left[ G\right] }\models \varphi $%
\textsl{.\hfill }$\square $

Of course $\models $ implies $\models _{\Omega }$.\ But the converse is
trivially false: there exists $\Omega $-valid formulae which are undecidable
in $ZFC$, for instance $\func{Con}\left( ZFC\right) $. Indeed, if $\left(
V_{\alpha }\right) ^{V\left[ G\right] }\models ZFC$ then $\left( V_{\alpha
}\right) ^{V\left[ G\right] }$ is a model of $ZFC$ and $\left( V_{\alpha
}\right) ^{V\left[ G\right] }\models \func{Con}\left( ZFC\right) $ So, $%
ZFC\models _{\Omega }\func{Con}\left( ZFC\right) $, but of course
(G\"{o}del) $ZFC\nvDash \func{Con}\left( ZFC\right) $.

It must be emphasized that $\Omega $-validity doesn't satisfy the \emph{%
compacity} property: there exist theories $T$ and formulae $\varphi $ s.t.
we have $T\models _{\Omega }\varphi $ even if for every finite subset $%
S\subset T$ we have $S\nvDash _{\Omega }\varphi $.\footnote{%
See Bagaria et al. [2005].}

By construction, $\Omega $-validity $\models _{\Omega }$ is itself
forcing-invariant:\footnote{%
See Woodin [2004].}

\textbf{Theorem (}$ZFC+PCW$\textbf{).} \textsl{If }$V\models ``T\models
_{\Omega }\varphi "$\textsl{\ then }$V[G]\models ``T\models _{\Omega
}\varphi " $\textsl{\ for every generic extension of }$V$\textsl{.}$\hfill
\square $

Woodin investigated deeply this new ``strong logic''. In particular he was
able to show that, under suitable $LCAs$, $CH$ ``rigidifies'' $V$ at the $%
\Sigma _{1}^{2}$-level ($\Sigma _{1}$ formulae for $V_{\omega +2}$):

\textbf{Theorem (Woodin, 1984)}. \textsl{Under }$PCW_{meas}$\textsl{\ (there
exists a proper class of measurable Woodin cardinals) and }$CH$\textsl{, }$%
\Omega $\textsl{-logic is generically complete at the }$\Sigma _{1}^{2}$%
\textsl{-level: for every }$\varphi $\textsl{\ of complexity }$\Sigma
_{1}^{2}$\textsl{\ either }$ZFC+CH\models _{\Omega }\varphi $\textsl{\ or }$%
ZFC+CH\models _{\Omega }\lnot \varphi $\textsl{. All generic extensions }$M$%
\textsl{\ and }$N$\textsl{\ of }$V$\textsl{\ satisfying both }$CH$\textsl{\
are }$\Sigma _{1}^{2}$\textsl{\ elementary equivalent.\hfill }$\square $

The metamathematical meaning of this result of \emph{conditional generic
absoluteness} is that if a problem is expressed by a $\Sigma _{1}^{2}$%
-formula $\varphi $ then it is ``settled by $CH$'' and immunized against
forcing under appropriate $LCAs.$ But:

\textbf{Theorem (Abraham, Shelah)}.\textsl{\ This is false at the }$\Sigma
_{2}^{2}$\textsl{\ level. For every large cardinal hypothesis }$A$\textsl{\
there exist generic extensions }$M$\textsl{\ and }$N$\textsl{\ satisfying
both }$CH$\textsl{\ s.t. in }$M$\textsl{\ there exist a }$\Sigma _{2}^{2}$%
\textsl{-wellorder of }$\mathbb{R}$\textsl{\ while in }$N$\textsl{\ all the }$%
\Sigma _{2}^{2}$\textsl{-subsets of }$\mathbb{R}$\textsl{\ are Lebesgue
measurable.\hfill }$\square $

In a second step, Woodin interpreted the $\Omega $-validity $T\models
_{\Omega }\varphi $ as the semantic validity for an $\Omega $-logic whose
syntactic derivation $T\vdash _{\Omega }\varphi $ had to be defined. His
idea was to witness the $\Omega $-proofs by particular sets that, under $PCW$%
, generalize the projective sets and can be interpreted without ambiguity in
every generic extension. It is the most difficult part of his work, not only
at the technical level but also at the philosophical level. The definition
(under $PCW$) is the following:

\textbf{Definition (}$PCW$\textbf{)}.\ $T\vdash _{\Omega }\varphi $\textsl{\
iff there exists a }universally Baire\textsl{\ (}$UB$\textsl{) set }$%
A\subseteq \mathbb{R}$\textsl{\ s.t. for every }$A$\textsl{-closed countable
transitive model (ctm) }$M$\textsl{\ of }$T$\textsl{\ we have }$M\models
\varphi $\textsl{\ (in other words }$M\models ``T\models _{\Omega }\varphi "$%
\textsl{).\hfill }$\square $

$A\subseteq \mathbb{R}$ is $UB$ if for every continuous map $f:K\rightarrow 
\mathbb{R}$ with source $K$ compact Hausdorff, $f^{-1}\left( A\right) $ has the
Baire property (there exists an open set $U$ s.t. the symetric difference $%
f^{-1}\left( A\right) \Delta U$ is meager). If $A\subseteq \mathbb{R}$ is $UB$
it is interpreted canonically in every generic extension $V[G]$ as $%
A_{G}\subseteq \mathbb{R}^{V[G]}$. This is due to the fact that there exists a 
\emph{tree presentation} of $A$. One identifies $\mathbb{R}$ with $\omega
^{\omega }$ and one considers trees $T\subset \left( \omega \times \gamma
\right) ^{\omega }$ and the projections $p[T]$ on $\omega ^{\omega }$ of
their infinite branches: 
\[
p[T]=\left\{ x\in \omega ^{\omega }\;\left| \;\exists z\in \gamma ^{\omega }%
\text{ with }\left( x\mid _{n},z\mid _{n}\right) \in T,\forall n\in \omega
\right. \right\} . 
\]

\noindent $A\subseteq \mathbb{R}$ is $UB$ iff there exist trees $T$ and $S$
s.t. $p[T]=A$ and $p[S]=\omega ^{\omega }-A$ in every generic extension $%
V[G] $. $p[T]$ yields a canonical interpretation of $A$ in every generic
extension $V[G]$.\ A ctm $M$ is called $A$-closed if, for every ctm $%
N\supseteq M$, $A\cap N\in N$, in particular for every generic extension $%
V[G]$ and $N=M[G]$ we have $A\cap M[G]\in M[G]$. If $A$ is Borelian, every
ctm is always $A$-closed.\footnote{%
If $M$ is $A$-closed for every $A$ of $\Pi _{1}^{1}$-complexity, then $M$ is
well-founded.}\ But it is no longer the case for general $UB$ sets.

As far as, in the definition of $T\vdash _{\Omega }\varphi $, the class of
admissible models is \emph{restricted} to $A$-closed ctms, logic becomes
strengthened. Of course, $T\vdash \varphi $ implies $T\vdash _{\Omega
}\varphi $, but the converse is false for the same reasons as for $\vDash
_{\Omega }$. Indeed, $ZFC\vdash _{\Omega }\func{Con}\left( ZFC\right) $
because every suitable ctm provides a model of $ZFC$ and validates therefore 
$\func{Con}\left( ZFC\right) $.

More technically, what is really needed for the definition $T\vdash _{\Omega
}\varphi $ are $UB$ sets $A\subseteq \mathbb{R}$ sharing the following two
properties:

\begin{enumerate}
\item  $L\left( A,\mathbb{R}\right) \vDash AD^{+}$, where $AD^{+}$ is a
strengthening of the axiom of determinacy saying that not only all $%
A\subseteq \mathbb{R\simeq \omega }^{\omega }$ are determined, but also all the 
$\pi ^{-1}\left( A\right) $ for all maps $\pi :\lambda ^{\omega }\rightarrow 
\mathbb{\omega }^{\omega }$ with an ordinal $\lambda <\frak{c}^{+}$;

\item  every $A\subset \mathcal{P}\left( \mathbb{R}\right) \cap L\left( A,\mathbb{R%
}\right) $ is $UB$.\footnote{%
See Woodin [2000].}
\end{enumerate}

\noindent $PCW$ implies these two properties and is therefore a good
hypothesis.

It must be emphasized that this definition of $\Omega $-provability is very
original. As explain Joan Bagaria, Neus Castells and Paul Larson in their ``$%
\Omega $-logic primer'':

\begin{quotation}
\noindent ``The notion of $\Omega $-provability differs from the usual
notions of provability, e.g., in first-order logic, in that there is no
deductive calculus involved.\ In $\Omega $-logic, the same $UB$ set may
witness the $\Omega $-provability of different sentences.\ For instance, all
tautologies have the same \emph{proof} in $\Omega $-logic, namely $\emptyset 
$.\ In spite of this, it is possible to define a notion of height of proof
in $\Omega $-logic.''
\end{quotation}

\noindent As Patrick Dehornoy explained to me (private communication), in $%
\Omega $-logic a proof $\vdash _{\Omega }\varphi $ is a certificate of some
property of the formula $\varphi $. This witnessing is no longer a
derivation iterating syntactic rules but a $UB$ subset of $\mathbb{R}$. What is
common to classical and $\Omega $-logics is that a very ``small'' object
endowed with a precise internal structure warrants the validity of $\varphi $
in a lot of immensely large models.

Woodin proved that $\Omega $-logic is \emph{sound}: if $T\vdash _{\Omega
}\varphi $ then $T\models _{\Omega }\varphi $, i.e. (under $PCW$) if $\vdash
_{\Omega }\varphi $ then $\models \varphi $ in all $ZFC$-models $\left(
V_{\alpha }\right) ^{V[G]}$. He then formulated the main conjecture:

$\Omega $-\textbf{conjecture (1999)}. $\Omega $\textsl{-logic is complete:
if }$\models _{\Omega }\varphi $\textsl{\ then }$\vdash _{\Omega }\varphi $%
\textsl{.\hfill }$\square $

As he emphasized in Woodin (2002, p.~517):

\begin{quotation}
\noindent ``If the $\Omega $-conjecture is true, then generic absoluteness
is equivalent to absoluteness in $\Omega $-logic and this in turn has
significant metamathematical implications''.
\end{quotation}

\noindent Indeed (Dehornoy, 2007), the $\Omega $-conjecture means that any
formula $\varphi $ valid in a lot of immensely large models satisfying $LCAs$
are certified by $UB$ subsets of $\mathbb{R}$. The key fact proved by Woodin is
the link of the concept of $\Omega $-derivability with the existence of 
\emph{canonical models} for $LCAs$ (that is models which are in a certain
way minimal and universal, as $L$ for $ZFC+CH$). The $\Omega $-conjecture
expresses essentially the hypothesis that every $LCA$ admits a canonical
model.

\textbf{Theorem.}\ $\vdash _{\Omega }\varphi $ \textsl{iff }$ZFC+A\vdash
\varphi $\textsl{\ for every large cardinal axiom }$A$\textsl{\ admitting a
canonical model.} $\hfill \square $

Now, the key point is that when $H_{2}$ is rigidified, $CH$ becomes
automatically \emph{false}.

\begin{quotation}
\noindent ``If the theory of the structure $\left\langle \mathcal{P}\left(
\omega _{1}\right) ,\omega _{1},+,.,\in \right\rangle $ is to be resolved on
the basis of a good axiom then necessarily $CH$ is false.''
\end{quotation}

\noindent The idea is that if the theory $T$ of $\mathcal{P}\left( \omega
_{1}\right) $ is completely unambiguous in the sense that there exists an
axiom $A$ s.t. $T\models \varphi $ iff $A\models _{\Omega }``T\models
\varphi "$, then $CH$ is necessarily false since the theory of $\mathcal{P}%
\left( \mathbb{R}\right) $ cannot share this property.

\textbf{Woodin theorem (}$2000$\textbf{, under }$PCW$\textbf{).} \textsl{(i)
For every ``solution'' for }$H_{2}$\textsl{\ (that is axioms freezing the
properties of }$H_{2}$\textsl{\ w.r.t. forcing) based on an }$\Omega $%
\textsl{-complete axiom }$A$\textsl{\ (i.e. for every }$\varphi \in H_{2}$%
\textsl{, either }$ZFC+A\vdash _{\Omega }``\left( H_{2},\in \right) \models
\varphi "$\textsl{\ or }$ZFC+A\vdash _{\Omega }``\left( H_{2},\in \right)
\models \lnot \varphi "$\textsl{), }$CH$\textsl{\ is false.\ (ii) If the }$%
\Omega $\textsl{-conjecture is valid, every ``solution'' for }$H_{2}$\textsl{%
\ is based on an }$\Omega $\textsl{-complete axiom and therefore }$CH$%
\textsl{\ is false.}$\hfill \square $

The proof uses Tarski results on the impossibility of defining truth and is
quite interesting (Woodin 2001, p.~688). Let 
\[
\Gamma =\left\{ \ulcorner \varphi \urcorner :ZFC+A\vdash _{\Omega }``\left(
H_{2},\in \right) \models \varphi "\right\} 
\]

\noindent be the (extremely complicated) set of G\"{o}del numbers of the
sentences $\Omega $-valid in $H_{2}$. By hypothesis, $\Gamma $ is $\Omega $%
-recursive in the sense there exists a $UB$ set $B$ s.t. $\Gamma $ is
definable and recursive in $L\left( B,\mathbb{R}\right) $. Now, $PCW$ implies
that $\Gamma $ being $\Omega $-recursive, it is definable in $\left( H\left( 
\frak{c}^{+}\right) ,\in \right) $.\ If $CH$ would be valid, then $\frak{c}
=\omega _{1}$, $H\left( \frak{c}^{+}\right) =H_{2}$ and $\Gamma $ would be
definable in $H_{2}$, which would violate Tarski theorem.

It is in that sense Woodin (2001, p.~690) can claim:

\begin{quotation}
\noindent ``Thus, I now believe the Continuum Hypothesis is solvable, which
is a fundamental change in my view of set theory''.
\end{quotation}

\section{Conclusion}

Hugh Woodin has already proved a great part of the $\Omega $-conjecture. 

Other approaches to the continuum problem in the set theoretical framework
of $LCAs$ have been proposed.\ One of the most interesting alternative is
provided by Matthew Foreman's (2003) concept of \emph{generic large cardinal}
($GLC$) defined by elementary embeddings $j:V\prec M$ of $V$ in inner models 
$M$ not of $V$ itself but of generic extensions $V[G]$ of $V$ . Such generic 
$LCAs$ can support rather $CH$ than $\lnot CH$.\ 

But all these results show what are the difficulties met in elaborating a
``good'' set theoretical determination of the continuum. The old Kantian
opposition between ``conceptual'' (symbolic) and ``intuitive'', or, in
Feferman's terms, between ``determined'' and ``inherently vague'', remains
insuperable.\ They justify some sort of G\"{o}del's platonism comprising
additional axioms as some kind of ``physical hypotheses''. The nominalist
antiplatonist philosophy of mathematics criticizing such axioms (in
particular $LCAs$) as ontological naive beliefs must be reconsidered and
substituted for a \emph{``conditional'' platonism} in Woodin's sense, a
platonism which would be ``conditional'' to axioms which ``rigidify'' the
continuum and make its properties forcing-invariant.

In my 1991, 1992 and 1995 papers on the continuum problem, I introduced the
concept of ``\emph{transcendental platonism}''. Classical platonism is a
na\"{i}ve realist thesis on the ontological independence of mathematical
idealities, and as such is always dialectically opposed to anti-platonist
nominalism. Even to day, the debates concerning the status of mathematical
idealities remain trapped into the realist/nominalist dialectic.\footnote{%
See e.g. Maddy [2005] on ``na\"{i}ve realism'', ''robust realism'', ``thin
realism'', etc.} The main achievement of transcendentalism has been to
overcome this scholastic antinomy between realism and nominalism and to show
that mathematical and physical objectivity were neither ontological nor
subjective.\ Objectivity is always \emph{transcendentally constituted} and
therefore \emph{conditional}, relative to eidetico-constitutive rules.\ A
platonism defined in terms of objectivity and not ontology, is a
transcendental platonism immune to the classical aporias of metaphysical
transcendent platonism. As far as the question of the continuum is
concerned, the eidetico-constitutive rules are the axioms of set theory and
transcendental platonism means that the continuum problem can have a well
determined solution in a ``rigid'' universe where $\mathbb{R}$ become \emph{%
conditionally} generically absolute.\ I think that Woodin's conditional
platonism can therefore be considered as a transcendental platonism relative
to the continuum problem.


\bigskip

\begin{thebibliography}{99}
\bibitem{}  Bagaria, J., Castells, N., Larson, P., 2005.\ ``An $\Omega $
-logic Primer'', \emph{Set Theory}, CRM 2003-2004, Birkhauser (2006), 1-28.

\bibitem{}  Becker, H., 1992. ``Descriptive Set Theoretic Phenomena in
Analysis and Topology``, \emph{STC [1992]}, 1-25.

\bibitem{}  Bell, J., 2005. \emph{The Continuous and the Infinitesimal in
Mathematics and Philosophy}, Polimetrica,
http://www.polimetrica.com/categories/02cat.html.

\bibitem{}  Bellotti, L., 2005.\ ``Woodin on the Continuum Problem: an
overview and some objections'', \emph{Logic and Philosophy of Science}, III,
1.

\bibitem{}  Benacerraf, P., Putnam, H. (eds.), 1964. \emph{Philosophy of
Mathematics: Selected Readings}, Prentice Hall, Englewood Ciffs, New-Jersey.

\bibitem{}  Breysse, O., De Glas, M., 2007. ``A New Approach to the Concepts
of Boundary and Contact: Towards an Alternative to Mereotopology'', \emph{%
Fundamenta Informaticae}, 78, 2, 217-238.

\bibitem{}  Chihara, Ch., S., 1990. \emph{Constructibility and Mathematical
Existence}, Clarendon Press, Oxford.

\bibitem{}  Cohen, P., 1963-1964.\ ``The independence of the Continuum
Hypothesis'', I, \emph{PNAS}, 50 (1963) 1143-1148 and II, \emph{PNAS}, 51
(1964) 105-110.

\bibitem{}  Cohen, P.\ 1966.\ \emph{Set Theory and the Continuum Hypothesis}%
, Benjamin, New York.

\bibitem{}  Dehornoy, P., 1989. ``La d\'{e}termination projective
d'apr\`{e}s Martin, Steel et Woodin'', \emph{S\'{e}minaire Bourbaki}, \#710.

\bibitem{}  Dehornoy, P., 2003.\ ``Progr\`{e}s r\'{e}cents sur
l'hypoth\`{e}se du continu (d'apr\`{e}s Woodin)'', \emph{S\'{e}minaire
Bourbaki}, \#915.

\bibitem{}  Dehornoy, P., 2007.\ ``Au-del\`{a} du forcing: la notion de
v\'{e}rit\'{e} essentielle en th\'{e}orie des ensembles'', \emph{Logique,
dynamique et cognition}, (J.B. Joinet, ed.), Publications de la Sorbonne,
Paris, 147--169.

\bibitem{}  Easton, W.\ B., 1970.\ ``Powers of regular cardinals'', \emph{AML%
}, 1 (1970), 139-178 (abridged version of the 1964 Ph.D. thesis).

\bibitem{}  Farah, I., Larson, P.\ B., 2005.\ ``Absoluteness for Universally
Baire Sets and the Uncountable, I'', \emph{Quaderni di Matematica}, 17
(2006), 47-92.

\bibitem{}  Feferman S., 1989. ``Infinity in Mathematics: Is Cantor
Necessary?'', \emph{Philosophical Topics}, XVII, 2, 23-45.

\bibitem{}  Feng, Q., Magidor, M., Woodin, W.\ H., 1992.\ ``Universally
Baire Sets of Reals'', \emph{STC 1992}, 203-242.

\bibitem{}  Foreman, M., Magidor, M., Shelah, S., 1986. ``$0^{\#}$ and Some
Forcing Principles'', \emph{Journal of Symbolic Logic}, 51 (1986), 39-47.

\bibitem{}  Foreman, M., 1998.\ ``Generic large cardinals: new axioms for
mathematics?'', \emph{ICM}, Vol. II, Berlin.

\bibitem{}  Foreman, M., 2003.\ ``Has the Continuum Hypothesis been
settled?'', http://www.math.helsinki.fi/logic/LC2003/presentations.

\bibitem{}  G\"{o}del, K., 1938. ``The consistency of the axiom of choice
and the generalized continuum hypothesis'', \emph{PNAS}, 25 (1938), 556-557.

\bibitem{}  G\"{o}del, K., 1940. ``The consistency of the axiom of choice
and of the generalized continuum hypothesis with the axioms of set theory'', 
\emph{Annals of Math. Studies}, Princeton University Press, Princeton.

\bibitem{}  G\"{o}del K., 1947. ``What is Cantor's Continuum Problem'', 
\emph{American Mathematical Monthly}, 54 (1947), 515-545 (reprinted in \emph{%
Benacerraf-Putnam [1964]}, 470-485).

\bibitem{}  Grigorieff, S., 1976. ``D\'{e}termination des jeux bor\'{e}liens
d'apr\`{e}s Martin'', \emph{S\'{e}minaire Bourbaki}, \#478.

\bibitem{}  Husserl, E., 1900-1901. \emph{Logische Untersuchungen}, Max
Niemeyer, Halle (1913).

\bibitem{}  Jackson, S., 1989. ``$AD$ and the very fine structure of $L(\mathbb{%
\ R})$'', \emph{Bulletin of the American Mathematical Society}, 21, 1 (1989)
77-81.

\bibitem{}  Jech, T., 1978.\ \emph{Set Theory}, Academic Press, San Diego.

\bibitem{}  Kanamori, A., 1995. \emph{The Higher Infinite}, Perspectives in
Mathematical Logic, Springer-Verlag, Berlin-Heidelberg-New York.

\bibitem{}  Kant, I. \emph{Kants gesammelte Schriften}, Preussische Akademie
der Wissenschaften, Georg Reimer, Berlin, 1911.

\bibitem{}  Maddy, P., 1988. ``Believing the Axioms I, II'', \emph{The
Journal of Symbolic Logic}, 53, 2 (1988), 481-511~; 53, 3 (1988), 736-764.

\bibitem{}  Maddy, P., 2005.\ ``Mathematical existence'', \emph{The Bulletin
of Symbolic Logic}, 11, 3 (2005), 351-376.

\bibitem{}  Martin, Donald, 1975. ``Borel Determinacy'', \emph{Annals of
Mathematics}, 102 (1975), 363-371.

\bibitem{}  Martin, David, 2001.\ \emph{Characterizations of }$0^{\#}$,
Master Thesis, Carnegie Mellon University.

\bibitem{}  Martin, D., Steel, J., 1989. ``A Proof of Projective
Determinacy'', \emph{Journal of the American Mathematical Society}, 2, 1
(1989), 71-125.

\bibitem{}  Moschovakis, Y., 1980. \emph{Descriptive Set Theory},
North-Holland.

\bibitem{}  Panza, M., 1992. ``De la Continuit\'{e} comme Concept au Continu
comme Objet'', \emph{Le Labyrinthe du Continu}, (J.-M.\ Salanskis, H.\
Sinaceur, eds.), Springer, Paris, 16-30.

\bibitem{}  Peirce, C.\ S., 1960.\ \emph{Collected Papers}, (C. Hartshorne
and P. Weiss, eds.), The Belknap Press of Harvard University\ Press,
Cambridge, Mass.

\bibitem{}  Peirce, C. S., 1976. \emph{New elements of mathematics}, (C.
Eisele, ed.), Mouton Publishers, The Hague.

\bibitem{}  Petitot, J., 1979. ``Infinitesimale'', \emph{Enciclopedia Einaudi%
}, VII, 443-521, Einaudi, Torino.

\bibitem{}  Petitot J., 1989. ``Rappels sur l'Analyse non standard'', \emph{%
La math\'{e}matique non-standard}, (H. Barreau, J. Harthong, eds.), Editions
du CNRS, Paris, 187-209.

\bibitem{}  Petitot J., 1991. ``Id\'{e}alit\'{e}s math\'{e}matiques et
R\'{e}alit\'{e} objective. Approche transcendantale'', \emph{Hommage \`{a}
Jean-Toussaint Desanti}, (G. Granel ed.), 213-282, Editions TER, Mauvezin.

\bibitem{}  Petitot J., 1992. ``Continu et Objectivit\'{e}. La
bimodalit\'{e} objective du continu et le platonisme transcendantal'', \emph{%
Le Labyrinthe du Continu}, (J.-M. Salanskis, H. Sinaceur, eds.),
Springer-Verlag, Paris, 239-263.

\bibitem{}  Petitot J., 1994. ``Phenomenology of Perception, Qualitative
Physics and Sheaf Mereology'', \emph{Philosophy and the Cognitive Sciences},
(16th Wittgenstein Symposium, R. Casati, B. Smith, G. White, eds), Verlag
H\"{o}lder-Pichler-Tempsky, Vienna, 387-408.

\bibitem{}  Petitot J., 1995. ``Pour un platonisme transcendantal'', \emph{\
L'objectivit\'{e} math\'{e}matique. Platonisme et structures formelles}, (M.
Panza, J-M. Salanskis, eds), 147-178, Masson, Paris.

\bibitem{}  Schimmerling, E., 2001.\ ``The ABC's of Mice'', \emph{The
Bulletin of Symbolic Logic}, 7, 4 (2001), 485-503.

\bibitem{}  Solovay, R. M., 1971. ``Real-Valued Measurable Cardinals'', 
\emph{Axiomatic Set Theory}, (D. Scott ed.), Proceedings of Symposia in Pure
Mathematics, Vol. XIII, Providence, AMS, 397-428.

\bibitem{}  STC, 1992. \emph{Set Theory of the Continuum}, (H. Judah, W.
Just, H. Woodin, eds.), Springer, Berlin-New York.

\bibitem{}  Steel, J., 2000.\ ``Mathematics Needs New Axioms'',
http://math.berkeley.edu/\symbol{126}steel.

\bibitem{}  Steel, J., 2001.\ ``Inner Model Theory'', to appear in the \emph{%
Handbook of Set Theory}.

\bibitem{}  Steel, J., 2004.\ ``Generic absoluteness and the continuum
problem'', http://math.berkeley.edu/\symbol{126}steel/talks/Lectures.html.

\bibitem{}  Stern, J., 1976. ``Le probl\`{e}me des cardinaux singuliers
d'apr\`{e}s Jensen et Silver'', \emph{S\'{e}minaire Bourbaki}, \#494.

\bibitem{}  Stern, J., 1984. ``Le probl\`{e}me de la mesure'', \emph{%
S\'{e}minaire Bourbaki}, \#632.

\bibitem{}  Thom, R., 1992. ``L'Ant\'{e}riorit\'{e} ontologique du Continu
sur le Discret'', \emph{Le Labyrinthe du Continu}, (J.-M.\ Salanskis, H.\
Sinaceur, eds.), Springer, Paris, 137-143.

\bibitem{}  Weyl, H., 1918. \emph{Das Kontinuum. Kritishe Untersuchungen
\"{u}ber die Grundlager der Analysis}, Veit, Leipzig, (English translation,
Lanham, University Press of America, 1987).

\bibitem{}  Woodin, W.\ H., 1999.\ \emph{The Axiom of Determinacy, Forcing
Axioms, and the Nonstationary Ideal}, De Gruyter, Berlin.

\bibitem{}  Woodin, W. H., 2000. ``The Continuum Hypothesis'', \emph{Logic
Colloquium}, Paris.

\bibitem{}  Woodin, W.\ H., 2001a.\ The continuum hypothesis, I-II, \emph{%
Notices Amer.\ Math.\ Soc.}, 48, 6 (2001), 567-576 and 48, 7 (2001) 681-690.

\bibitem{}  Woodin, W. H., 2001b.\ ``The $\Omega $ Conjecture'', \emph{%
Aspects of Complexity}, De Gruyter, Berlin, 155-169.

\bibitem{}  Woodin, W.\ H., 2002.\ ``Beyond $\Sigma _{1}^{2}$
Absoluteness'', \emph{ICM 2002}, Vol. I, 515-524.

\bibitem{}  Woodin, W.\ H., 2003. ``Set theory after Russell: the journey
back to Eden'', \emph{One Hundred Years of Russell's Paradox}, (G.\ Link,
ed.), De Gruyter, Berlin, 29-47 (2004).

\bibitem{}  Woodin, W.\ H., 2004.\ ``Is there really any evidence that the
Continuum Hypothesis has no answer?'',
http://www.lps.uci.edu/home/conferences/
Laguna-Workshops/LagunaBeach2004/irv.pdf
\end{thebibliography}
\end{document}